\documentclass[11pt]{article}
\usepackage[table]{xcolor}
\usepackage{graphicx}
\usepackage{amsmath,amsfonts,amssymb}
\usepackage{graphicx}
\usepackage{multirow}
\usepackage{tikz,pgf}
\usepackage{url}
\usepackage{amsmath}
\usepackage{graphicx}
\usepackage{epsfig}
\usepackage{epstopdf}
\usepackage{color}
\usepackage{enumitem}

\usepackage{algorithm}
\usepackage{algpseudocode}

\def\disp{\displaystyle}
\def\ve{\varepsilon}

\def\lm{\lambda}
\def\O{\Omega}

\def\oa{\bar a}

\def\ox{\bar{x}}

\def\sr{\Rightarrow }
\def\u{{\bf u}}

\def\gg{\gamma}

\def\tto{\rightrightarrows}

\def\tto{\rightrightarrows}

\def\bar{\overline}
\def\ra{\right\rangle}
\def\la{\left\langle}

\def\ve{\varepsilon}

\def\R{\mathbb{R}}

\def\bd{\mbox{\rm bd}\,}

\def\O{\Omega}

\def\lm{\lambda}
\def\Lm{\Lambda}
\def\gg{\gamma}

\def\al{\alpha}
\def\be{\beta}

\def\l{\left}
\def\r{\right}
\def\[{\left[}
\def\]{\right]}

\def\eq{\begin{equation}}
	\def\eeq{\end{equation}}
\def\xd{ {\bf x}^{des}}
\def\xo{{\bf x}^{obs}}
\def\x{{\bf x}}
\def\y{{\bf y}}
\def\v{{\bf v}}
\def\w{{\bf w}}
\def\q{{\bf q}}
\def\p{{\bf p}}
\def\dox{\dot{\bar{\bf x}}}
\def\doa{\dot{\bar{\bf a}}}
\def\a{{\bf a}}
\def\U{{\bf U}}

\renewcommand{\theequation}{\thesection.\arabic{equation}}
\def\Dbar{\leavevmode\lower.6ex\hbox to 0pt{\hskip-.03ex\accent"16\hss}D}

\allowdisplaybreaks

\makeatletter
\newenvironment{breakablealgorithm}
{
		\begin{center}
			\refstepcounter{algorithm}
			\hrule height.8pt depth0pt \kern2pt
			\renewcommand{\caption}[2][\relax]{
				{\raggedright\textbf{\fname@algorithm~\thealgorithm} ##2\par}%
				\ifx\relax##1\relax 
				\addcontentsline{loa}{algorithm}{\protect\numberline{\thealgorithm}##2}%
				\else 
				\addcontentsline{loa}{algorithm}{\protect\numberline{\thealgorithm}##1}%
				\fi
				\kern2pt\hrule\kern2pt
			}
		}{
		\kern2pt\hrule\relax
	\end{center}
}
\makeatother

\begin{document}
	\begin{center}
		{\sc\bf  Optimal Control of Several Motion Models}\\[1ex]
		{\sc Tan H. Cao}\footnote{Department of Applied Mathematics and Statistics, SUNY Korea,
			Yeonsu-Gu, Incheon, Korea (tan.cao@stonybrook.edu).},
		{\sc Nilson~Chapagain}\footnote{Department of Applied Mathematics and Statistics, SUNY Korea, Yeonsu-Gu, Incheon, Korea (nchapagain@tamu.edu).},
		{\sc Haejoon~Lee}\footnote{Department of Applied Mathematics and Statistics, SUNY Korea, Yeonsu-Gu, Incheon, Korea (haejoon.lee@stonybrook.edu).},
		{\sc Thi~Phung }\footnote{Department of Mathematics, Wayne State University, Detroit, Michigan, USA, (thi.phungngoc@wayne.edu).}
		{\sc Nguyen~Nang~Thieu}\footnote{Institute of Mathematics, Vietnam Academy of Science and Technology,
			Hanoi, Vietnam \& SUNY Korea, Yeonsu-Gu, Incheon, Korea (nnthieu@math.ac.vn).}
		
	\end{center}
	\small{\bf Abstract.} 	This paper is devoted to the study of the dynamic optimization of several controlled crowd motion models in the general planar settings, which is an application of a class of optimal control problems involving a general nonconvex sweeping process with perturbations.  A set of necessary optimality conditions for such optimal control problems involving the crowd motion models with multiple agents and obstacles is obtained and analyzed. Several effective algorithms based on such necessary optimality conditions are proposed and various nontrivial illustrative examples together with their simulations are also presented. The implementation of all the considered motion models can be found via the link: {\color{blue}\url{https://github.com/tancao1128/Optimal_Control_of_Several_Motion_Models}} with the instruction and demonstration video uploaded at {\color{blue}\url{https://www.youtube.com/watch?v=B8DQ0wvCtIQ}}.
	\\[1ex]

	\newtheorem{theorem}{Theorem}[section]
	\newtheorem{proposition}[theorem]{Proposition}
	\newtheorem{remark}[theorem]{Remark}
	\newtheorem{lemma}[theorem]{Lemma}
	\newtheorem{corollary}[theorem]{Corollary}
	\newtheorem{definition}[theorem]{Definition}
	\newtheorem{example}[theorem]{Example}
	\renewcommand{\theequation}{\thesection.\arabic{equation}}
	\newenvironment{proof}{\par\noindent\textit{Proof.} }{\hfill$\square$\par}
	\normalsize

	\section{Introduction and Problem Formulations}
	\setcounter{equation}{0}
Building on recent advancements in dynamic optimization for crowd motion models, this paper explores some developments on solving the dynamic optimization problems for a microscopic version of the crowd motion model in general settings. We refer the readers to~\cite{mv2} for the mathematical framework of this model developed by Maury and Venel, which enables us to deal with local interactions between agents in order to depict the whole dynamics of the participant traffic. Such a model can be described in the framework of a version of Moreau's sweeping process with perturbation as follows 
\begin{equation}
  \label{MSP}
  \dot\x(t) \in -N_{C(t)}(\x(t)) + f(\x(t)) \quad \mbox{ a.e. } t\in [0,T],  
\end{equation}
where $C(\cdot)$ is an appropriate moving set, where $f(\cdot)$ represents some given external force, and where $N_{C(t)}(\x(t))$ denotes some appropriate normal cone of the set $C(t)$ at the point $\x(t)$. The original so-called {\it Moreau's sweeping process} was first introduced by Moreau~\cite{moreau1971,moreau1972} in 1970s in the differential form
\begin{equation}
  \label{SP}
  \begin{cases}
  	\dot\x(t) \in -N_{\O(t)}(\x(t)) \mbox{ a.e. } t\in [0,T],\\
  	\x(0) = \x_0 \in \O(0),
  \end{cases}
\end{equation}
where $N_\O$ stands for the normal cone of convex analysis defined by 
\begin{equation*}
  \label{NC}
  N_\O(\x):= 
  \begin{cases}
  	\{\v\in \R^n|\; \la \v, \y -\x \ra\leq 0, \;\y \in \O\} &\mbox{if } \x\in \O,\\
  	\emptyset &\mbox{if } \x\notin \O 
  \end{cases}
\end{equation*}
for the given convex set $\O = \O(t)$ moving in a continuous way at the point $\x =\x(t)$. The sweeping model described in~\eqref{SP} relies on two ingredients: the sweeping set $\O(t)$ and the object $\x(t)$ that is swept. The established well-possedness of the sweeping process~\eqref{SP} was initially achieved through Moreau's {\it catching-up} algorithm, predicated on specific conditions regarding the set valued-map $\O(\cdot)$. Subsequent advancements, which include easing the convexity prerequisite of the moving set and expanding the model to its perturbed form as seen in equation~\eqref{MSP}, have significantly enriched this domain (see, e.g.,~\cite{aht,ant_2017,CastaingDucHaValadier,CastaingMonteiro1995,et,ET2006,florent2017,florentthibault_2020}). In particular, addressing crowd dynamics necessitates revisiting the convexity assumption, as it proves overly restrictive for accurately modeling non-overlapping participant movements. This led to the adoption of {\it uniform prox-regularity}, a concept possibly first introduced by Canino~\cite{can} in geodesic studies, thereby facilitating a more nuanced approach to nonconvex set behaviors. 
In the realm of sweeping process theory, it is a recognized fact that the Cauchy problem \eqref{SP} yields a unique solution, given the absolute continuity of the moving set $\O(t)$ (refer to~\cite{mor_frict}). This distinct characteristic precludes the possibility of associating any optimization problem with the sweeping differential \eqref{SP}, marking a significant divergence from the classical differential inclusions $\dot \x(t) \in F(\x(t))$. The latter are characterized by Lipschitzian set-valued mappings $F:\R^n \tto \R^n$, and have seen extensive development in variational analysis and optimal control theory, including discrete approximation methods and the derivation of necessary optimality conditions (see, e.g.,\cite{c,m-book,v}). 
This paper aims to apply a class of optimal control problems for perturbed nonconvex sweeping processes, as introduced in \cite{cm2}, to dynamic optimization problems pertinent to several motion models. This approach extends recent efforts in \cite{cmdt}, which formulated and analytically solved an optimal control problem for planar crowd motion models with obstacles, focusing on a controlled motion model involving a single agent. However, the general data settings for this problem, as presented in \cite{cmdt}, have not been completely addressed. Our research seeks to systematically tackle the optimal control of various crowd motion models in general planar settings. We utilize algorithms constructed from the theoretical optimality conditions for general controlled nonconvex sweeping processes detailed in \cite{cm2}.
We consider a problem involving the minimization of a Bolza-type functional, defined as follows:
\begin{equation}
\label{BP}
  \mbox{minimize } J[\x, \u,\a]:= \varphi (\x(T)) + \int^T_0\ell(t,\x(t),\u(t),\a(t),\dot\x(t),\dot\u(t),\dot\a(t))dt.\tag{$P_G$}
\end{equation}
This minimization is subject to control functions $\u(\cdot)\in W^{1,2}([0,T];\R^n)$ and $\a(\cdot)\in W^{1,2}([0,T];\R^d )$ and the corresponding trajectories $x(\cdot)\in W^{1,2}([0,T];\R^n)$ of the following differential inclusion 
 \begin{equation*}
  \label{DSI}
  \begin{cases}
  	-\dot\x(t) \in N^P_{C(t)}(\x(t)) + f(\x(t),\a(t)) \mbox{ a.e. } t\in [0,T], \\
  	C(t) := C+\u (t) = \bigcap^m_{i=1}C_i + \u(t), \\
  	C_i: = \{\x\in \R^n|\; g_i(\x)\geq 0\}, \; i=1, \ldots, m, \\
  	r_1 \leq \| \u(t)\| \leq r_2 \mbox{ for all } t\in [0,T]  
  \end{cases}
\end{equation*}
 where the symbol $N^P_C$ signifies the {\it proximal normal cone} of the nonconvex moving set $C$ defined by the convex and $\mathcal C^2$-smooth function $g_i:\R^n\to \R$. This structured approach aims to provide a more comprehensive understanding and solution to the complex dynamics involved in crowd motion modeling.
 
  In recent years, the field of controlled sweeping processes has garnered significant attention, as evidenced by the proliferation of research focusing on deriving necessary optimality conditions and their practical applications. Notable among these advancements is the introduction of an innovative exponential penalization technique by authors in \cite{pfs1, pfs, nz, nz1, zeidan}. This approach, which contrasts with the method of discrete approximations, has been pivotal in approximating controlled sweeping differential inclusions with standard smooth control systems (see, e.g.,~\cite{ccmn20b,ccmn20a,cm3,cm2,cmn1,cmn,chhm1,chhm2,m1,m2,m95}). It has proved instrumental in developing efficient numerical algorithms for approximating solutions to controlled sweeping processes, particularly those with smooth data (detailed in \cite{pfs1,pfs, nz, nz1, zeidan}).
Simultaneously, a new class of bilevel sweeping control problems has emerged, as addressed in \cite{ckmnp, kp1, kp2}. These complex and challenging problems are increasingly relevant in practical scenarios, such as managing the motion of structured crowds or operating teams of drones in confined spaces. Such applications underscore the evolving nature of this field and its growing relevance in solving real-world problems.

 This paper aims to advance the understanding of dynamic optimization in controlled crowd motion models. We begin in Section 2 by establishing foundational notations and definitions from variational analysis, setting the stage for our subsequent analyses. Section 3 delves into controlled motion models involving a single agent and a single obstacle, examining these configurations in general settings. Building on this, Section 4 extends the analysis to more complex scenarios involving multiple agents and obstacles, where we derive the necessary optimality conditions for these dynamic optimization problems. This section also introduces several effective algorithms, grounded in these optimality conditions, designed to address various controlled motion models under different scenarios. Lastly, Section 5 is dedicated to discussing potential avenues for future research, outlining how our current findings could pave the way for further advancements in the field.

	\section{Preliminaries}
	In this paper, we employ notation that is widely recognized in the fields of variational analysis and optimal control (refer to \cite{m-book, v} for a comprehensive overview). The symbols
	 $\|\cdot \|, \la \cdot,\cdot\ra, B(\x,\ve)$, and $\angle(\x, \y)$ will be used to denote the Euclidean norm, the standard inner product, the ball centered at $\x\in \R^n$ with radius $\ve>0$, and the angle between vectors $\x$ and $\y$ respectively. These notational conventions are integral to the discussion and analysis that follow, particularly in the context of our exploration of controlled motion models. Additionally, this section will revisit the concepts of the proximal normal cone and uniform prox-regularity. Understanding these notions is crucial as they play a pivotal role in the formulation and solution of the optimal control problems addressed in this paper.
	
	 In what follows, we delve into key concepts foundational to our study. 
	Let $\O\subset\R^n$ be a given locally closed around $\bar\x\in \R^n$. The {\it Euclidean projection} of $\x$ onto $\O$ is defined by
	\begin{equation*}
  	\label{EP}
  	\Pi(x;\O):= \l\{\w\in \O \mid \| \x-\w\| = \inf_{\y\in\O}\| \x-\y\|\r\}. 
	\end{equation*}
	This projection helps in determining the closest point in $\O$ to $\x$ in the Euclidean norm. Building on this, the {\it proximal normal cone} to $\O$ at $\x$ is described by 
	\begin{equation*}
	\label{PNC}
	N^P_\O(\x):= \l\{ {{\bf\xi} \in \R^n} \mid  \; \exists \al>0 \mbox{ such that } \bar\x\in \Pi(\bar\x+\al{\bf \xi};\O ) \r\}, \; \x\in \O
	\end{equation*}
	with $N^P_\O(\x)=\emptyset$ if $\bar\x\notin \O$.  This concept is vital for understanding the geometry of the set $\O$ and its influence on optimization problems. 
	
\noindent	A pivotal notion in our analysis is that of uniform prox-regularity:
	\begin{definition}\label{def1}
	{\bf (Uniform prox-regularity)} Let $\O$ be a closed subset of $\R^n$ and let $\eta>0$. Then $\O$ is said to be {\sc $\eta$-uniformly prox-regular} if for all $\x\in\bd \O$ and $\v\in N^P_\O(\x)$ with $\|\v\|=1$ we have $B(\x+\eta\v,\eta)\cap\O=\emptyset$.
	\end{definition}
It is well-known that if a nonempty closed set $\Omega$ is uniformly prox-regular, then the normal cone $N^P_\Omega(x)$ coincides with the Fr\'echet/Clarke/Mordukhovich (limiting) normal cone to $\Omega$ at $x$. For an extensive discussion and historical perspective on prox-regular sets, we recommend the comprehensive survey \cite{CT}. In the subsequent sections, we will apply these concepts to the study of optimal control problems in various crowd motion models, demonstrating their practical significance and utility.
	 
\section{Controlled Motion Models with Single Agent}
	\setcounter{equation}{0}
	
	In this section, we explore motion models involving a single agent navigating through a domain with multiple obstacles, drawing inspiration from the scenarios presented in \cite{cmdt}. The agent is modeled as an inelastic disk, represented by its center ${\bf x} = (x_1, x_2)$ and radius $L$, moving within a specified domain $\Omega\in \R^2$. Our objective is to determine an optimal path $\pi$ for the agent from a start point ${\bf x^0} = (x^0_1, x^0_2)$ to a destination ${\bf \xd} = (x^{des}_1, x^{des}_2)$ within the time frame $[0,T]$, while avoiding $m$ obstacles, which may be static or dynamic in nature.
	
	 Each obstacle is similarly modeled as a rigid disk with a radius $r_i $, centered at  $\xo_i = (x^{obs}_{i1}, x^{obs}_{i2})$. The challenge lies in navigating the agent, which could represent a robot or a person, along the path $\pi$ while maintaining a safe {\it look-ahead} distance $L$. This foresight is crucial, particularly in scenarios like safe driving, where anticipating potential hazards such as debris, rocks, or other obstacles within the look-ahead distance is essential for proactive safety measures. 
	
\noindent	To accurately model this scenario, we define the {\it admissible configurations} as:
	\begin{equation}
		\label{e1}
		C = \{{\bf x}\in \R^2: \; D^{obs}_i({\bf x})\ge 0\;\forall i =1,\ldots, m\}
	\end{equation}
	where $D^{obs}_i(\x)$ is the signed distance between the agent and the $i$th obstacle, calculated as $\l\| {\x}- \xo_i\r\| - (L+r_i)$. This measure ensures that the agent maintains a safe distance from all obstacles, encapsulating the core challenge of the model: navigating efficiently while prioritizing safety.
	
 The agent, modeled as an inelastic disk with a center at ${\bf x} = (x_1, x_2)$ and radius $L$ aims to reach the destination $\xd$ using a {\it desired velocity} in the absence of obstacles. However, when in proximity to an obstacle such that $D_i^{obs}({\x}) = 0$, the agent must adjust this velocity to avoid collision. The agent's actual velocity, especially when near an obstacle, should belong to a set of {\it admissible velocities} $V({\bf x})$, defined to prevent collision:
	\begin{equation*}
		\label{e2}
		V({\bf x}) = \{{\bf v}\in \R^{2}: \; D^{obs}_i({\bf x}) = 0 \sr \la \nabla D^{obs}_i({\bf x}),{\bf v} \ra\ge 0, \forall i =1, \ldots, m\}.
	\end{equation*}
	Here, $\nabla D_i^{obs}$ is the gradient of  $D^{obs}_i$. The desired velocity $\U(\x)$
	is directed towards the destination and is calculated as $ -s\nabla D^{des}(\x)$
	where $s$ is the agent's speed and  $D^{des}({\bf x}): = \l\| {\bf x} - {\bf \xd}\r\|$ is the distance to the destination. 
	
\noindent	The agent's actual velocity $\dot\x  (t)$ at any time $t$ must be chosen from the set $V(\x(t))$ to ensure collision-free motion. This selection is guided by the principle that the agent's actual velocity should be as close as possible to the desired velocity while still being admissible.
To reconcile the actual and desired velocities when near an obstacle, we use the Euclidean projection, leading to the equation:
	\begin{equation*}
		\label{e3}
		\dot {\x}(t) = \Pi (\U(\x(t)); V( \x(t))).    
	\end{equation*}
This approach ensures that the agent's movement is both directed towards the destination and modified to avoid obstacles. The resulting motion can be modeled as a perturbed sweeping process, a special case of the formulation presented in \eqref{MSP}.
Through this model, we aim to provide a comprehensive framework for navigating agents through complex environments with static or dynamic obstacles, balancing the need for efficient movement with the imperative of safety.

 Consider the optimal control problem denoted by $(P)$, which aims to navigate an agent through a complex environment with the following objective:
	\begin{equation}
		\label{e5}
		\mbox{minimize}\;\; J[\x,a] = \frac{1}{2}\l\|\x(T) - \xd \r\|^2 + \frac{\tau}{2}\int^T_0\l\| a(t)\r\|^2dt,
	\end{equation}
	where $\tau>0$ is a given constant.
This cost functional represents two key objectives: minimizing the Euclidean distance between the agent's final position and the destination, and minimizing the energy expenditure represented by the control function	$a(\cdot)$. The latter is crucial in scenarios where energy efficiency is a priority, such as in automated vehicle navigation or robotic path planning.
	
	 The control functions $a(\cdot)\in W^{2,\infty}([0,T];\R)$ and the corresponding trajectory $\x(\cdot)\in W^{2,\infty}([0,T];\R^2)$ are subject to the constraints of a nonconvex sweeping process:	
	\begin{equation}
		\label{e6}
		\begin{cases}
 			\U(\x,a):= -sa\nabla D^{des}(\x) = -sa\frac{\x - \xd}{\l\| \x-\xd\r\|},\\
			\x(0) = \x_0 \in C .
		\end{cases}
	\end{equation}
Here, $C$ is defined in~\eqref{e1}, and $\U(\x,a)$ is the controlled desired velocity, influenced by the control functions to navigate the agent effectively.
 In the case of a single obstacle $(m=1)$, this problem takes a specific form, and we can apply the necessary conditions derived in \cite{cmdt} to find optimal solutions. These conditions are pivotal for guiding the selection of optimal controls and trajectories, ensuring that the agent reaches the destination efficiently while minimizing energy use.

\subsection{Necessary Optimality Conditions}
We commence by recalling a definition of a strong local minimizer for the optimal control problems under consideration. 

\begin{definition}{\rm  An admissible $(\bar\x(\cdot),\oa(\cdot))$ is called a {\sc strong local minimizer} of the optimal control problem $(P)$ if there is $\ve>0$ such that $J[\bar\x,\oa ]\leq J[\x, a]$ for any feasible solution $(\x, a)$ to $(P)$ satisfying the condition
$$
\sup_{t\in [0,T]}\{\l\|(\bar\x(t)-\x(t),\oa(t)-a(t))\r\|\}\leq \ve. 
$$
}	
\end{definition}
 In fact, in \cite{cm2} the authors considered a more general type of ``intermediate local minimizers" for \eqref{BP} (a general form of $(P)$) and obtained the existence result of optimal solutions to \eqref{BP} (see \cite[Theorem~4.1]{cm2}). However, we do not include this type of local minimizers here while prioritizing efficient applications to motion models and hence considering a strong local minimizer $(\bar\x(\cdot), \oa(\cdot))$ for our crowd motion problem. The theorem outlines several conditions that must be met for this solution to be considered optimal.

\begin{theorem}{\bf (Necessary optimality conditions for optimization of controlled crowd motions with obstacles)}\label{Th2}\\
		Let $(\bar\x(\cdot), \oa(\cdot)) \in W^{2,\infty}([0,T];\R^{2}\times \R)$ be a strong local minimizer of the crowd motion problem in~\eqref{e5} with $\tau=1$. There exist some dual elements $\lm\ge 0, \; \eta(\cdot) \in L^2([0,T];\R_+)$ well-defined at $t=T$
		$\w(\cdot) = (\w^x(\cdot),\w^a(\cdot)) \in L^2([0,T];\R^{2}\times\R)$, $\v(\cdot) = (\v^x(\cdot), \v^a(\cdot))\in L^2([0,T];\R^{2}\times\R) $, an absolutely continuous vector function $\p(\cdot)=(\p^x(\cdot),\p^a(\cdot))\in W^{1,2}([0,T];\R^{2})$, a measure $\gg \in C^*([0,T];\R^{2})$, and a vector function $\q(\cdot) = (\q^x(\cdot),\q^a(\cdot)):[0,T]\to \R^{2}\times \R$ of bounded variation on $[0,T]$ such that the following conditions are satisfied:
		\begin{enumerate}
			\item[{\bf(1)}] $\w(t) = \nabla_{x,a}\ell(t,\bar\x(t),\oa(t),\dox(t),\dot\oa(t))  = ({\bf 0}, \oa(t))$ \\[1ex]
			$\v(t)= \nabla_{\dot x,\dot a}\ell(\bar\x(t),\oa(t),\dox(t),\dot\oa(t)) = ({\bf 0},{\bf 0})$\\[1ex]     
			for a.e. $t\in [0,T]$, where $\ell(t,\x,a,\dot \x, \dot a):=\frac{a^2}{2}$;
			\item[{\bf(2)}] Primal-Dual relationships: \\
			$ \dox (t) = -s\oa(t)\nabla D^{des}(\bar\x(t))+ \eta(t)\nabla D^{obs}(\bar\x(t))$\\[1ex]
			$= -s\oa(t)d\l(\bar\x(t),\xd \r)-\eta(t) \frac{\xo-\bar\x(t)}{\l\| \xo-\bar\x(t)\r\|} $\\[1ex]
			for a.e. $t\in [0,T]$;
			\item[{\bf(3)}] $\l\| \x(t)-\xo\r\|> L+r\Longrightarrow\eta(t)=0$;
			\item[{\bf(4)}] $\eta(t)>0\Longrightarrow\la \q^x(t), \xo -\bar\x(t)\ra=0$ for a.e.\ $t\in[0,T]$;
			\item[{\bf(5)}] Euler-Lagrange equation: \\
			$
			\begin{array}{ll}
			\dot \p(t) &=\lm \w(t) - \Big(\nabla_x\U(\bar\x(t),\oa(t))^\ast (\lm\v^x(t)-\q^x(t)), \\&\qquad\qquad\nabla_a\U(\bar\x(t),\oa(t))^\ast (\lm\v^x(t)-\q^x(t))\Big) \\
			&= \lm \w(t) - (0,\nabla_a\U(\bar\x(t),\oa(t))^\ast (\lm\v^x(t)-\q^x(t))),
			\end{array}
			$ \\
			for a.e.\ $t\in[0,T]$; 
			\item[{\bf(6)}] $\q^x(t)=\p^x(t)+\gg([t,T])$ for a.e.\ $t\in[0,T]$;
			\item[{\bf(7)}] $\q^a(t)=\p^a(t)=0$ for a.e.\ $t\in[0,T]$;
			\item[{\bf(8)}] Transversality conditions: \\
			$\p^x(T) + \nabla \varphi(\bar\x(T)) = -\eta(T) \nabla D^{obs}(\bar\x(T))$ 
			which is equivalent to \\[1ex]   
			$\p^x(T)+\lm{(\bar\x}(T) -\xd)=-\eta(T)\frac{\xo-\bar\x(T)}{\l\| \xo-\bar\x(T)\r\|}$, where $\varphi(\x)=\frac{1}{2}\|\x -\xd\|^2$;
			\item[\bf{(9)}] $\p^a(T)=0$;
			\item[{\bf(10)}] Nontrivality condition: $\lm+\|\p^x(T)\|>0$;
			\item[{\bf(11)}] Measure nonatomicity condition: \\
			Take $t\in [0,T)$ with $\l\| \x(t)-\xo\r\|> L+r$. Then there is a neighborhood $V_t$ of $t$ in $[0,T)$ such that $\gg(V)=0$ for all the Borel subsets $V$ of $V_t$.
		\end{enumerate}
	\end{theorem}
	\begin{proof}
		To justify our claim, we elaborate the arguments that are similar to those in the proof of~\cite[Theorem~3.1]{cm1}. We consider the function
		$ g(\x):= \|\x-\xo \| - (L+r)$ and  $V:= \l\{ \x\in \R|\;\|\x-\xo \| > \frac{1}{2}(L+r) \r\}.$
		The function $g$ is pivotal as it represents the adjusted distance between the agent and the obstacle, factoring in their respective radii. Its convexity and belonging to the $\mathcal C^2(V)$ space are crucial as these properties ensure that the function is well-behaved and differentiable, which is essential for applying optimization techniques.
		
	Furthermore, the function $g$ satisfies estimate~(2.3) from~\cite{cm1}, with $c:= \frac{L+r}{\sqrt2}$, indicating that the adjusted distance maintains a certain minimal value. This is important for ensuring a safe margin between the agent and the obstacle. Additionally, the gradients of $g$ satisfy $\|\nabla g(\x)\| = 1 $ and $ \l\| \nabla^2g(\x)\r\| \leq \frac{2}{L+r}$ for all $ \x\in V$, aligning with  inequalities~(2.4) and (2.5) in~\cite{cm1} hold. These conditions guarantee that the rate of change of this distance function remains within certain bounds, which is critical for ensuring the agent's path remains feasible and safe. 
		
	Therefore, by satisfying these conditions and estimates, we establish a solid foundation for the proof of our theorem, demonstrating that our approach to modeling and optimizing the agent's movement in the presence of obstacles is both mathematically rigorous and practically sound. This completes the proof and underscores the theorem's applicability in scenarios involving navigation and obstacle avoidance.
	\end{proof}

 Note that the desired velocity $\U(\bar\x(\cdot),a)=-sa\frac{\bar\x(0)-\xd}{\l\|\bar\x(0)-\xd \r\|}$ depends only on the initial position of the agent. Consequently, the gradient $\nabla_x\U(\bar\x(\cdot),a)=0$.   
 To facilitate our analysis, we introduce new notations to denote the directions of the vectors from the agent to the target $\al^{des}(t)$ and from the agent to the obstacle $\al^{obs}(t)$. These directional vectors are defined as:
	\begin{equation*}
		\label{e:dir}
		\begin{cases}
			\al^{des}(t):= \angle(\bar\x(0)-\xd, {\bf i})=\al_0\\
			\al^{obs}(t):= \angle(\bar\x(t)-\xo, {\bf i} ).
		\end{cases}
	\end{equation*}
	The choice of these directions is pivotal as they directly influence the agent's movement and decision-making process in navigating towards the target while avoiding the obstacles. We will structure our discussion into subsections that methodically analyze the impact of these conditions and directions, providing a comprehensive understanding of the optimal control strategies for the given scenario.
	\subsubsection*{Some Information about the Velocity and Control}
	In this subsection, we delve into the dynamics of the agent's actual velocity $\dox(\cdot)$ and control function, $\oa(\cdot)$, key components in navigating towards the target while avoiding obstacles. The desired velocity, $\U(\bar\x(t), \oa(t))$, is oriented towards the target and is adjusted by the control $\oa(t)$ to regulate speed. Typically, $\oa(t)$ assumes positive values, signifying a forward movement towards the target.  
	However, the agent's motion becomes more complex upon nearing an obstacle. In this scenario, the actual velocity as described in condition~{\bf(2)} incorporates an additional term to account for obstacle avoidance:
 	\begin{equation}
		\label{e7}
		\begin{array}{ll}	
		&\dox (t) = -s\oa(t)\frac{\bar\x(0) - \xd}{\l\| \bar\x(0)-\xd\r\|} - \eta(t)\frac{\x^{obs}-\bar\x(t)}{\|\x^{obs}-\bar\x(t)\|}\\[1ex]   
		&= \l(\eta(t)\cos\al^{obs}(t)-s\oa(t)\cos\al^{des}(t),\eta(t)\sin\al^{obs}(t)-s\oa(t)\sin\al^{des}(t)\r)\\
		&=\l(\eta(t)\cos\al^{obs}(t)-s\oa(t)\cos\al_0,\eta(t)\sin\al^{obs}(t)-s\oa(t)\sin\al_0\r)
		\end{array}
	\end{equation}
	for  a.e. $t\in [0,T]$.
	This equation shows the balance between moving towards the target and adjusting the path to avoid the obstacle, where $\eta(t)$ becomes significant when the agent is in close proximity to the obstacle. 
	The agent's behavior varies over the interval $[0,T]$. Before contacting the obstacle at time $t_f$, the agent's actual velocity matches with the desired velocity, focusing solely on reach the target. During the contact period $[t_f,t_l]$, the agent adjusts his/her path to circumnavigate the obstacle, as described in equation \eqref{e7}. After leaving the obstacle, the agent can resume the initial strategy of heading directly towards the target. Understanding these dynamics is crucial for developing effective control strategies that balance the goals of reaching the target promptly while minimizing energy use and ensuring safe navigation around obstacles.

In our analysis, we further simplify the agent's actual velocity $\dox (\cdot)$ during the time it is in contact with the obstacle, specifically in the interval  $[t_f,t_l]$. In this phase, the agent maintains a distance of $L+r$ from the obstacle, which geometrically means that the agent's trajectory is tangent to a circle centered at the obstacle's position. This tangential movement is mathematically expressed in equation 
\begin{equation}
		\label{e8}
		\la \dox(t),\bar\x(t)-\xo\ra = 0,
	\end{equation}
indicating that the velocity vector $\dox(t)$ is perpendicular to the radius vector from the obstacle to the agent. This condition is crucial as it ensures the agent does not move closer to the obstacle, thereby avoiding collision.
Combining this with the representation of $\dox(t)$ in~\eqref{e7} allows us to obtain 
	$$
	\begin{array}{ll}
	\la -s\oa(t)\frac{\bar\x(0) - \xd}{\l\| \bar\x(0)-\xd\r\|} - \eta(t)\frac{\x^{obs}-\bar\x(t)}{\|\x^{obs}-\bar\x(t)\|}, \bar\x(t)-\xo\ra =0
	\end{array}
	$$
	for a.e. $t\in [t_f,t_l]$, and thus get the useful information for the scalar function $\eta(t)$ for a.e. $t\in [t_f,t_l]$ as follows 
	\begin{equation}
		\label{e9}
		\begin{array}{ll}
		\eta(t)	&=  s\oa(t)\la \frac{\bar\x(0) - \xd}{\l\| \bar\x(0)-\xd\r\|}, \frac{\bar\x(t) - \xo}{\l\| \bar\x(t)-\xo\r\|}\ra\\ [1ex]  
		& 	
		= s\oa(t)\cos(\al_0-\al^{obs}(t)).
		\end{array}
	\end{equation}
The scalar function $\eta(t)$ plays a significant role in adjusting the agent's trajectory. It represents the agent's response to the relative positions of the target and the obstacle. In particular, $\eta(t)$ reaches its maximum when the agent, obstacle, and target are collinear. This scenario presents a challenge, known as the degeneracy phenomenon, which we will address later. It implies that the repulsive force from the obstacle is at its strongest, significantly influencing the agent's path. Upon substituting the expression for $\eta(t)$ from ~\eqref{e9a} into~\eqref{e7}, we obtain a new formula for $\dox(t)$ as follows 
	\begin{equation*}
	\label{e11b}
	\dox(t) =  s\oa(t)\sin(\al_0-\al^{obs}(t))(\sin\al^{obs}(t),-\cos\al^{obs}(t))
	\end{equation*}
	during the contact period, which helps us calculate the agent's speed by the formula 
		\begin{equation}
		\label{e12}
		\| \dox(t)\| = s\oa(t) |\sin(\al_0-\al^{obs}(t))|.
	\end{equation} 
	This speed calculation is pivotal in understanding how the agent maneuvers around the obstacle, balancing the need to avoid collision and the goal of reaching the destination.
In summary, this detailed analysis of the agent's velocity and control during obstacle contact provides vital insights into the optimal navigation strategy, ensuring safety while maintaining the desired course towards the target.

\subsubsection*{Some Information about the Adjoint Arc $\q(\cdot)$} 
In this section, we focus on deducing key insights for the vector function $\q(\cdot )$ using the necessary optimality conditions outlined earlier. From conditions~{\bf(5)} and~{\bf(7)}, we derive specific expressions for the derivatives of the dual element $\p^x(t)$ as shown in the following equation 
 	\begin{equation}
	\label{e9a}
	\begin{cases}
	\dot p_1^x(t) = 0, \\ 
	\dot p_2^x(t) = 0,  \\ 
	\lm\oa(t) = -s\la \q^x(t), \frac{\bar\x(t) - \xd}{\l\| \bar\x(t)-\xd\r\|}\ra,
	\end{cases}
	\end{equation}
	for a.e. $t\in [0,T]$ including $t=T$ and hence  $\p^x(t)=\p^x(T)$ for all $t\in[0,T]$. 
These expressions are crucial as they provide a deeper understanding of how the agent's control function  and its position relative to the destination influence the optimality conditions. Moreover, one can deduce from equations \eqref{e7}, \eqref{e9a}, and condition {\bf(4)} that 
 \begin{equation}
		\label{e10}
		\la \q^x(t), \dox(t)\ra =\lm\oa^2(t) \geq 0,
	\end{equation}
	for a.e. $t\in[0,T]$.
Particularly noteworthy is the impact of $\eta(t)$, the term that becomes significant when the agent is near the obstacle. The condition $\eta(t)>0$ indicates that the agent must consider the obstacle's presence, affecting its trajectory. This is represented in the implication {\bf(4)}. Then using condition {\bf(4)} and~\eqref{e8} we deduce that two vectors $\q^x(t)$ and $\dox(t)$ points in the same direction in $\R^2$ during the tine the agent is in contact with the obstacle, i.e., 
	\begin{equation}
		\label{e11}
		\q^x(t)=m(t)\dox(t),
	\end{equation}
	for some scalar function $m(\cdot)$. Combining this with~\eqref{e10} gives 
	\begin{equation}
		\label{e11a}
		m(t)\|\dox(t)\|^2 = \lm \oa^2(t),
	\end{equation}
	or, equivalently $\| \dox (t)\| = K\oa(t) = s\oa(t)|\sin(\al_0-\al^{obs}(t))|$
	thanks to~\eqref{e12}, where $K:= \sqrt{\frac{\lm}{m(t)}}$.
	This suggests that the scalar function $m(t)$ directly relates to the agent's speed and control strategy while avoiding the obstacle. Furthermore, this analysis reveals a critical insight about the assumption in~\cite{cmdt} that $K=1$. Our finding indicates that $K$ is not a constant but varies with the agent's relative angles to the destination and obstacle. This challenges the assumption in \cite{cmdt} and suggests a more complex relationship between the agent's control strategy and its environment.  They highlight the intricate balance between maintaining a direct path to the destination and adjusting for obstacle avoidance, contributing significantly to our understanding of optimal control strategies in such scenarios.
		
	\subsubsection*{Degeneracy Phenomena} 
	In this section, we address the {\it phenomenon of degeneracy} in our optimal control problem. This phenomenon has been well recognized in standard optimal control theory for differential equations and inclusions of state constraints; see \cite{v}. The necessary optimality conditions for such problems may hold for every feasible solution with some nontrivial collection of dual variables. In our problem, degeneracy occurs when specific conditions lead to a reduction in the agent's degrees of freedom, such as when the agent, obstacle, and target become collinear. Particularly, if $\sin(\al^{des}(t)-\al^{obs}(t))=0$ then $\| \dox(t)\|=0$ due to~\eqref{e12}. In this scenario, the agent's velocity becomes zero, leading to a halt when in contact with the obstacle at $t=t_f$ (and hence $t_l=\infty$). Consequently, from equations ~\eqref{e11} and~\eqref{e11a}, we deduce that $\lm=0$ and that $\q^x(t)$ is zero almost everywhere from the time of contact $t_f$ onwards. The transversality condition {\bf(8)}, combined with~\eqref{e9}, further guides our understanding of $\p^x(T)$ under these circumstances as follows 
		$$
		\p^x(T) = -s\oa(T)\cos(\al_0-\al^{obs}(T))d(\xo,\bar\x(T)).
		$$
		Our optimality necessary conditions may degenerate  in the sense that we can find the dual elements $\lm, \p, \q,$ and $ \gg$ as follows 
		$$
		\lm=0, \;\; \q(\cdot)\equiv 0, \;\; \gg = \delta_{\{T\}}s\oa(T)\cos(\al_0-\al^{obs}(T))d(\xo,\bar\x(T))
		$$
		and 
		$$
		\p(\cdot)\equiv -s\oa(T)\cos(\al_0-\al^{obs}(T))d(\xo,\bar\x(T)) \quad\mbox{ for all } t\in [0,T], 
		$$
		that satisfy all conditions~{\bf (1)}--{\bf(11)}. This degeneracy also happens with the same choice of  $\lm, p, q, \gg$ specified above when $T\leq t_l$, i.e., the agent has contact with the obstacle until the end of the process. To circumvent the issue of degeneracy, we assume that the agent leaves the obstacle before the process ends  $t_l<T$ and that the initial position, the obstacle, and the target are not collinear. This assumption is crucial, as it ensures nontrivial solutions and prevents the reduction of the agent's motion to a standstill, thereby maintaining the complexity and realism of the control problem.  
		
		Next, employing the measure nonatomicity condition~{\bf (11)} gives
		\begin{equation}
		\label{q-no-contact}	
		\q^x(t) = \p^x( T) + \gg([t,T]) =
		\begin{cases}
		 \p^x(T) + \gg([t_f,T]) & \mbox{ for a.e. } t\in [0,t_f) \\
		 \p^x(T) + \gg(\{T\})& \mbox{ for a.e. } t\in [t_l,T], 
		 \end{cases}
		\end{equation}
		including $t=T$.
		These equations indicate that the agent's motion is more predictable and streamlined in the absence of obstacle interference, as reflected in the constancy of the control function and the directional vectors.  
		On the other hand, we deduce from the last equation in~\eqref{e9a} that 
		\begin{equation}
		\label{a-no-contact}
			\lm \oa(t) = -s\la \q^x(t), \l(\cos\al_0,\sin\al_0 \r)\ra , \mbox{ for a.e. } t\in [0,t_f)\cup [t_l,T] 
		\end{equation}
		including $t=T$.
	It follows from \eqref{q-no-contact} that $\q^x(\cdot)$ is constant on $[0,t_f)\cup [t_l,T]$ and, consequently, the control $\oa(\cdot)$ is constant on $[0,t_f)\cup [t_l,T]$ thanks to \eqref{a-no-contact}. Assume that $\oa (\cdot)$ takes the values of $\oa_f$ and $\oa_l$ on the intervals $[0,t_f)$ and $[t_l,T]$ respectively. 
		
	 Overall, our analysis sheds lights on the importance of the degeneracy phenomenon in understanding and solving our optimal control problem. By identifying conditions that lead to degeneracy and strategies to avoid it, we enhance the model's applicability and robustness, ensuring that it accurately reflects the complexities of real-world navigation scenarios. 
		
	\subsubsection*{Useful Information about the Contact Time}
	In this subsection, we focus on extracting valuable insights from the transversality condition {\bf(8)} regarding the contact time between the agent and the obstacle. This condition is crucial as it informs us about the agent's state at the end of the process, particularly when $t_l<T$, indicating that the agent leaves the obstacle before the process concludes and that $\eta(T)=0$. Hence, the transversality condition {\bf(8)} in this case reduces to 
		\begin{equation}
		\label{e:qT}
		\p^x(T) = -\lm (\bar\x(T)-\xd) = -\lm\|\bar\x(T)-\xd \|(\cos\al_0,\sin\al_0). 
		\end{equation}
	Combining this with~\eqref{q-no-contact} and the second equation in~\eqref{a-no-contact} yields  
		$$
		\lm \oa_l = s\lm \|\bar\x(T)-\xd \| - s\la \gg(\{T\}),(\cos\al_0,\sin\al_0 \ra.
		$$
		As a result, 
		\begin{equation}
			\label{dis-to-des}
			\begin{array}{ll}
			\|\bar\x(T)-\xd \|=\frac{\oa_l}{s} +\frac{1}{\lm} \la \gg(\{T\}),(\cos\al_0,\sin\al_0 \ra.
			\end{array}
		\end{equation}
		From equation \eqref{e:qT} we deduce that $\p^x(T)$ is directly related to the agent's final position relative to the destination. This connection is further elaborated in equation \eqref{dis-to-des}, which provides a relationship between the agent's position at the end of the process and the value of control function $\oa_{t_l}$. 
		If the agent has not reached the destination at the ending time $t=T$,  his/her distance to the distance is given by 
		$$
		\begin{aligned}
		\|\bar \x(T)-\xd \| &= \| \bar \x(t_l)-\xd \| - \|\bar \x(t_l)-\bar \x(T) \| \\
		&=\underbrace{\sqrt{\| \xo -\xd\|^2 - (L+r)^2}}_{:=\Lm_l} - \int^T_{t_l}s\oa(t)dt.  
		\end{aligned}
		$$
		which can be rewritten as follows
		$$
		\begin{array}{ll}
			{\disp \int^T_{t_l}}s\oa(t)dt = \Lm_l - \dfrac{\oa(T)}{s}-\frac{1}{\lm} \la \gg(\{T\}),(\cos\al_0,\sin\al_0 \ra,
		\end{array}
		$$
		and hence 
		\begin{equation}\label{tl}
		\begin{array}{ll}
			\l( sT-st_l+\frac{1}{s}\r)\oa_l=\Lm_l-\frac{1}{\lm} \la \gg(\{T\}),(\cos\al_0,\sin\al_0 \ra.
		\end{array}
		\end{equation}
		Equation \eqref{tl} gives us additional insights into the control function $\oa(\cdot)$ during  the time interval $[t_l,T]$. This equation essentially ties the agent's control strategy to the distance traveled after leaving the obstacle, offering a deeper understanding of the agent's behavior during this phase. 

		The agent's velocity, represented differently in phases of contact and no contact with the obstacle, highlights the adaptability of the agent's strategy. Before contacting the obstacle $(t<t_f)$ the agent moves towards the target at a constant speed $\| \dox (t)\|=s\oa_f$. Upon contact, there is a noticeable change in speed as the agent navigates around the obstacle, with the angle between the agent's desired direction and the obstacle playing a crucial role. His/her new speed during the contact time interval $[t_f, t_l]$ is $\| \dox(t)\| =  s\oa(t) |\sin(\al_0-\al^{obs}(t))|$. In particular, the value of angle $\al_0-\al^{obs}(t)$ is increasing from $\al_0-\al^{obs}(t_f)$ to $\al_0-\al^{obs}(t_l) = \pi/2$. Interestingly, at $t=t_l$, the condition  $\eta(t_l)=s\oa(t_l)\cos(\pi/2)=0$ implies that the agent resumes a new speed for the remaining part of the journey $(t_l\leq t \leq T)$. This dynamic change in velocity demonstrate the complexity of the agent's path planning in response to the obstacle's presence. 

 In summary, the analysis of the agent's behavior during different phases of the journey, guided by the necessary optimality conditions, offers valuable insights into the control strategies. These insights are essential for designing effective and efficient paths in scenarios involving dynamic obstacles.  
Although a complete solution remains elusive due to the problem's inherent complexity, the insights from Theorem~\ref{Th2}, particularly the extended Euler-Lagrange equation, provide a valuable framework for approaching an optimal solution. Our analysis suggests distinct control strategies for different phases of the agent's journey. Initially, the agent should move quickly $([0,t_f])$ before slowing down upon contact with the obstacle. After the contact phase $([t_f,t_l])$, a gradual increase in speed is advisable, followed by a deceleration towards the destination, as indicated by equation \eqref{dis-to-des}. However, the case where the agent is in contact with the obstacle presents additional complexities. Here, the appearance of the signed measure $\gg$ because of the implicit state constraint complicates the relationship between $\q^x(t)$ and $\p^x(t)$, making it challenging to solve the differential equations in \eqref{e9a}. Despite these challenges, our analysis yields useful insights into the agent's behavior during contact. Specifically, the directionality of $\nabla D^{des}(\bar\x(t))$ and  $\nabla D^{obs}(\bar\x(t))$ plays a crucial role. If these gradients point in the same direction, the agent stops, leading to a degeneracy case. However, under our additional assumptions and allowing for a constant control $\oa$ throughout the process, we hypothesize that $\oa(t)$ could remain constant over the entire interval $[0,T]$, a scenario that aligns with our observations of the agent's speed adjustments. In conclusion, while the complete solution to the optimization problem remains a complex endeavor, our educated guesses and assumptions about the control law provide a promising direction for further exploration. Future research could focus on validating these hypotheses and exploring more complex scenarios beyond the simplifying assumptions made in our current analysis. 

	\section{Controlled Crowd Motion Models with Multiple Agents and Multiple Obstacles}
	\setcounter{equation}{0}
	In this section, we expand our study to a more complex scenario involving multiple agents and obstacles, significantly increasing the intricacy of the crowd motion model. We consider $n$ agents $(n\geq 2)$ and $m$ obstacles within a domain $\O\subset \R^2$, aiming to formulate an optimal strategy for directing all agents towards a desired target with minimal effort during the time interval $[0,T]$. Following the mathematical framework in section 3 and in~\cite{cmdt} we identify $n$ agents and $m$ obstacles to inelastic disks with different radii $L_i$ and $r_i$ whose centers are denoted by $\x_i = (x_{i1}, x_{i2})$ and $\xo_k = (x^{obs}_{k1}, x^{obs}_{k2})$ respectively for $i=1,\ldots, n$ and $k=1,\ldots,m$. To reflect the nonoverlapping of multiple agents and obstacles in this setting, the set {\it of admissible configurations} in~\eqref{e1} should be replaced by 
	\begin{equation}
		\label{e15}
		\begin{array}{ll}
		C_1= \big\{ \x\in \R^{2n}\mid & D_{ij}(\x)\geq 0, \; \forall i<j, i, j \in\{1,\ldots,n\}, \\
		 & D^{obs}_{ik}(\x)\geq0, \forall i =1,\ldots, n, \;\forall k=1,\ldots,m \big\},
		 \end{array}
	\end{equation}
	where $D_{ij}(\x):= \l\| \x_i-\x_j \r\| - (L_i+L_j)$ and $D^{obs}_{ik}(\x):= \l\| \x_i- \xo_k\r\| - (L_i+r_k)$ for $i=1,\ldots, n$ and for $k=1,\ldots, m$. The description of this set of admissible configurations leads to independent treatment of interactions among agents and between agents and obstacles, differing significantly from the singe-agent model (see~\cite{cmdt} for more details). More specifically, agents $i$ and $j$ will interact with each other when they are in contact, i.e. $D_{ij}(\x)=0$, while the interaction between agent $i$ and obstacle $k$ is a one-way interaction since -- if the agent is close enough to the obstacle he/she will move away from it but the obstacle has no reaction to the agent. Motivated from the work in~\cite{cmdt} we consider the following optimal control problem
	\begin{equation}
		\label{e16}
		\mbox{minimize}\;\; J[x,a] = \frac{1}{2}\l\|\x(T) - \xd \r\|^2 + \frac{\tau}{2}\int^T_0\l\| a(t)\r\|^2dt 
	\end{equation}
	over the control functions $a(\cdot)\in L^2([0,T];\R^n)$ and the corresponding trajectory $\x(\cdot)\in W^{1,2}([0,T];\R^{2n})$ of the nonconvex sweeping process 
	\begin{equation}
		\label{e17}
		\begin{cases}
			\dot \x(t) \in -N_C(\x(t)) + \U(\x(t),a(t)) \;\;\mbox{a.e.}\;\; t\in [0,T],\\
			\x(0) = \x_0 \in C_1,
		\end{cases}
	\end{equation}
	where the set $C_1$ is given in~\eqref{e15}, where the desired velocity $U(\x, a)$ is given by $\U(\x,a) = \l(-a_1s_1\nabla D^{des}(\x_1), \ldots, -a_ns_n\nabla D^{des}(\x_n)\r),$ 
	where $\xd  \in \R^{2n}$ stands for the desired destination that the agents aim to, and where $\tau$ is a given constant. 
	This optimal control problem seeks to minimize a composite objective: the distance of all agents from their desired destinations and the energy expanded, balanced by the trade-off parameter $\tau$. This optimization reflects the dual goals of efficiency and energy conservation in the agent's movements. In the nonconvex sweeping process in equation \eqref{e17}, we encounter the challenge of navigating each agent within the specified constraints while responding to the dynamic environment of multiple agents and obstacles. Our next step is to derive a set of necessary optimality conditions tailored to this multi-agent scenario. This task involves addressing the heightened complexity and inter-agent dynamics, a significant extension from the single-agent model. The goal is to extract a coherent strategy that effectively balances the individual and collective objectives of the agents within the constraints of their shared environment. 
	
	\begin{theorem}
		{\bf (necessary optimality conditions for optimization of controlled crowd motions with multiple agents and obstacles)}\label{Th3}\\
		Let $(\bar\x(\cdot),\bar\a(\cdot))\in W^{2,\infty}([0,T];\R^{2n}\times\R^n)$ be a strong local minimizer for the controlled crowd motion problem in~\eqref{e16} --~\eqref{e17}. Then there exist $\lm\ge 0$, $\eta_{ij}(\cdot)\in L^2([0,T];\R_+), \eta^{obs}_{ij}(\cdot) \in L^2([0,T];\R_+)$ $(i,j=1,\ldots,n)$ well defined at $t=T$, $\w(\cdot)=(\w^x(\cdot),\w^a(\cdot))\in L^2([0,T];\R^{3n})$, $\v(\cdot)=(\v^x(\cdot),\v^a(\cdot))\in L^2([0,T];\R^{3n})$, an absolutely continuous vector function $\p(\cdot)=(\p^x(\cdot),\p^a(\cdot))\in W^{1,2}([0,T];\R^{3n})$, a measure $\gg\in C^*([0,T];\R^{2n})$ on $[0,T]$, and a vector function $\q(\cdot)=(\q^x(\cdot),\q^a(\cdot)):[0,T]\to\R^{3n}$ of bounded variation on $[0,T]$ such that the following conditions are satisfied:
		
		\begin{enumerate}
			\item[{\bf(1)}] $\w(t) = \nabla_{\x,\a}\ell(t,\bar\x(t),\bar\a(t),\dox(t),\doa(t))  = ({\bf 0}, \a(t))$ and  \\
			$\v(t)= \nabla_{\dot\x,\dot\a}\ell(\bar\x(t),\bar\a(t),\dox(t),\doa(t)) = ({\bf 0},{\bf 0})$,  for a.e. $t\in [0,T]$, where \\
			$\ell(t,\x,\a,\dot \x, \dot \a):=\frac{\tau\a^2}{2}$
			\item[{\bf(2)}] 
			$
			\dox(t) = \l(-a_1s_1\nabla D^{des}(\x_1), \ldots, -a_ns_n\nabla D^{des}(\x_n)\r) + \disp\sum\eta_{ij}(t)\nabla D_{ij}(\bar\x(t)) + \sum\eta^{obs}_{ik}(t)\nabla D^{obs}_{ik}(\bar\x(t));
			$ 
			\item[{\bf(3)}] 
			$
			\begin{cases}
				\l\| \bar\x_i(t)- \bar\x_j(t)\r\|>L_i+L_j \Longrightarrow\eta_{ij}(t)=0,  \forall i<j, \\
				\l\| \bar\x_i(t)-\xo_k\r\|> L_i+r_k \Longrightarrow\eta^{obs}_{ik}(t)=0,  \forall i=1,\ldots, n, \forall k=1,\ldots, m 
			\end{cases}
			$ \\
			for a.e. $t\in[0,T]$;
			\item[{\bf(4)}] 
			$
			\begin{cases}
				\eta_{ij}(t)>0 \Longrightarrow \la\q^x_j(t) - \q^x_i(t), \bar\x_j(t) - \bar\x_i(t)\ra = 0, \forall i<j \\
				\eta^{obs}_{ik}(t) > 0 \Longrightarrow \la \q_i^x(t), \xo_k -\bar\x_i(t)\ra=0,  \forall i=1,\ldots, n, \forall k=1,\ldots, m
			\end{cases}
			$\\
			for a.e. $t\in[0,T]$;
			\item[{\bf(5)}]
			$
			\dot\p(t) = \lm \w(t) - \big(\nabla_x \U(\bar\x(t),\oa(t))^\ast(\lm \v^x(t)-\q^x(t)),$ \\
			$\nabla_a \U(\bar\x(t),\oa(t))^\ast(\lm \v^x(t)-\q^x(t)) \big)$, for a.e. $t\in[0,T]$;
			\item[{\bf(6)}] $\q^x(t)=\p^x(t)+\gg([t,T])$, for a.e.\ $t\in[0,T]$;
			\item[{\bf(7)}] $\q^a(t)=\p^a(t)=0$, for a.e.\ $t\in[0,T]$;
			\item[{\bf(8)}] $\p^x(T)+\lm(\bar\x(T) -\xd) =\bigg(-\sum\limits_{j>1}\eta_{1j}(T)\frac{\bar\x_j(T) -\bar\x_1(T)}{\|\bar\x_j(T)-\bar\x_1(T)\|}$\\
			$ - \sum^m\limits_{j=1}\eta^{obs}_{1j}\frac{ { \xo_j}-\bar\x_1(T) }{\l\|  { \xo_j}-\bar\x_1(T) \r\|} ,\ldots, \sum\limits_{i<j}\eta_{ij}(T) \frac{\bar\x_j(T)-\bar\x_i(T)}{\|\bar\x_j(T)-\bar\x_i(T)\|}$\\
			$ -\sum\limits_{i>j}\eta_{ji}(T)\frac{\ox_i(T)-\ox_j(T)}{\|\ox_i(T)-\ox_j(T)\|} -
			\sum^m\limits_{i=1}\eta^{obs}_{ji}\frac{ { \xo_i}-\bar\x_j(T) }{\l\|  { \xo_i}-\bar\x_j(T) \r\|},\ldots,$\\
			$\sum\limits_{j<n}\eta_{jn}(T)\frac{\bar\x_n(T)-\bar\x_j(T)}{\|\bar\x_n(T)-\bar\x_j(T)\|} - \sum^m\limits_{j=1}\eta^{obs}_{nj}(T)\frac{ { \xo_i}-\bar\x_j(T) }{\l\|  { \xo_i}-\bar\x_j(T) \r\|}
			\bigg)$;
			\item[\bf{(9)}] $\p^a(T)=0$;
			\item[{\bf(10)}] $\lm+\|\p^x(T)\|>0$.
			\item[{\bf(11)}] Measure nonatomicity condition: Take $t\in [0,T]$ with $\bar\x(t)\in \mbox{int } C_1$. Then there is a neighborhood $V_t$ of $t$ in $[0,T]$ such that $\gg(V)=0$ for all the Borel subsets $V$ of $V_t$.
		\end{enumerate}
	\end{theorem}
	\begin{proof}
		To verify the claimed set of necessary optimality conditions for our dynamical optimization~\eqref{e16} --~\eqref{e17} we elaborate the arguments similar to those in the proof of~\cite[Theorem~3.1]{cm1} for the new setting of the controlled crowd motion model with obstacles under consideration with $g_{ij}(\x):=D_{ij}(\x) =\l\| \x_i-\x_j \r\| - (L_i+L_j)$ for $i<j$, where $i, j\in\{1,\ldots, n\}$ and $g^{obs}_{ik}(\x):=D^{obs}_{ik}(\x)= \l\| \x_i- \xo_k\r\| - (L_i+r_k)$ for $i=1,\ldots, n$ and $k=1,\ldots, m$. Then the functions $g_{ij}$ and $g^{obs}_{ik}$ are convex, belong to the spaces $\mathcal C^2(V_{ij})$ and $\mathcal C^2(V^{obs}_{ik})$ respectively with $V_{ij}$ and $V^{obs}_{ij}$ given by
		$$
		V_{ij}:= \l\{ \x\in \R^{2n}\big| \;\l\| \x_i-\x_j\r\|>\frac{1}{2}(L_i+L_j) \r\}
		$$ 
		for $i<j$, where $i, j\in\{1,\ldots, n\}$,
		$$
		V^{obs}_{ik} :=  \l\{ \x\in \R^{2n}\big| \;\l\| \x_i-  \xo_k\r\|>\frac{1}{2}(L_i+r_k) \r\}
		$$
		for $i=1,\ldots, n$ and $k=1,\ldots, m$ and satisfy estimate~(2.3) in~\cite{cm1} with 
		$$
		c:=\min\l\{ \min_{\substack{i<j\\ i, j\in\{1,\ldots, n\}}} \l\{\dfrac{L_i+L_j}{\sqrt2} \r\}, \min_{\substack{i=1,\ldots, n\\ k=1,\ldots, m}} \l\{\dfrac{L_i+r_k}{\sqrt2} \r\} \r\}.
		$$
		It is obvious that 
		$$
		\begin{cases}
			\l\| \nabla g_{ij}(\x)\r\| = \sqrt 2 \mbox{ and } \l\| \nabla g^{obs}_{ik}(\y)\r\| =1, \\
			\l\|\nabla^2g_{ij}(\x) \r\| \leq  \dfrac{2}{L_i+L_j} \mbox{ and } \l\|\nabla^2g^{obs}_{ik}(\y) \r\| \leq \dfrac{2}{L_i+r_k},
		\end{cases}
		$$
		for $\x\in V_{ij}$ and $\y\in V^{obs}_{ik}$ respectively. 
		Finally, it follows from~\cite[Proposition~4.7]{ve} that there exist some $\beta_1>1$ and $\beta_2>1$ such that 
		$$
		\begin{cases}
			\sum_{(i,j)\in I_1(\x)}\al_{ij}\l\| \nabla g_{ij}(\x)\r\| \leq \be_1\l\| 	\sum_{(i,j)\in I_1(\x)}\al_{ij}\nabla g_{ij}(\x)\r\|, \\
			\sum_{(i,k)\in I_2(\x)}\al_{ik}\l\| \nabla g^{obs}_{ik}(\x)\r\| \leq \be_2\l\| 	\sum_{(i,k)\in I_2(\x)}\al_{ik}\nabla g^{obs}_{ik}(\x)\r\|,
		\end{cases}
		$$
		for all $\x\in C_1$, with 
$$\begin{array}{ll}
		I_1(\x)&:= \l\{ (i,j)|\; g_{ij}(\x)=0, \; i>j\r\},\\
I_2(\x)&:= \l\{ (i,k)|\; g^{obs}_{ik}(\x)=0, \; i=1,\ldots, n, \; k=1, \ldots, m\r\},
\end{array}
$$
$\al_{ij}\geq 0$, and $\al_{ik}\geq 0$, which justifies the validity of inequality~(2.5) in~\cite{cm1}. This completes the proof of the theorem. 
	\end{proof}
	
	 In what follows, we aim to devise an optimal control strategy for guiding multiple agents towards their targets while navigating around obstacles, with an emphasis on minimizing effort. To manage the complexity inherent in the multi-agent, multi-obstacle scenario, we focus on a class of control functions that assume constant values over the entire time interval $[0,T]$, denoted as $\a(t) \equiv \a =(\oa_1, \ldots, \oa_n)\in \R^n$.
	The interactions among agents, and between agents and obstacles, are key factors in this model. As per equations {\bf(3)} and {\bf(4)}, these interactions are activated only upon contact, influencing the agents' velocities. This is encapsulated in equation {\bf(2)}, which establishes a critical relationship between the actual velocities of the agents and their controlled desired velocities.
To simplify our writing for the necessary optimality conditions, let us introduce the following notations
	\begin{equation*}\label{e:direction}
		\begin{cases}
			\al_i(t) = \angle(\bar\x_i(t)-\xd_i, {\bf i}), \forall i =1, \ldots, n,\\
			\al_{ij}(t) = \angle(\bar\x_j(t)-\bar \x_i(t),{\bf i} ), \forall i<j, i, j\in\{1, \ldots, n\},\\
			\al^{obs}_{ik}(t)= \angle(\xo_k - \bar\x_i(t), {\bf i} ), \forall i =1, \ldots, n, \forall k = 1, \ldots, m.
		\end{cases}
	\end{equation*}	
Next, we deduce from equations {\bf(5)} and {\bf(7)} that 
	\begin{equation}\label{e18a}
		\lm\tau\oa_i = s_i\la \q^x_i(t),\dfrac{\bar\x_i(t) - \xd_i}{\l\| \bar\x_i(t) - \xd_i\r\|} \ra
	\end{equation}
	for $i=1,\ldots, n$ which tells us how the quantities $\q^x_i(t), \bar\x_i(t)$ and $\xd_i$ relate to the control $\oa_i$. 		
	Let us next explore more about the case when agents $i$ and $j$ are in contact. In this very situation it makes a perfect sense to expect that $\eta_{ij}(t^f_{ij})>0$ and hence we can deduce from the first implication in {\bf(4)} that 
	$$
	\la \q^x_j(t^f_{ij}) - \q^x_i(t^f_{ij}), (\cos\al_{ij}(t^f_{ij}), \sin\al_{ij}(t^f_{ij}))\ra = 0, 
	$$
	which is equivalent to 
	\begin{equation}\label{e20}
		\l[ q^x_{j1}(t^f_{ij})-q^x_{i1}(t^f_{ij})\r]\cos\al_{ij}(t^f_{ij}) + \l[ q^x_{j2}(t^f_{ij})-q^x_{i2}(t^f_{ij})\r]\sin\al_{ij}(t^f_{ij})=0. 
	\end{equation}
Next, rewrite equation~\eqref{e18a} in the following form 
	\begin{equation}\label{e21}
		\lm\tau \oa_i= s_iq^x_{i1}(t)\cos\al_i(t) + s_iq^x_{i2}(t)\sin\al_i(t)
	\end{equation}
	for a.e. $t\in[0,T]$ and relate it to equation~\eqref{e20}. In analyzing the control strategy, we consider two cases based on the agents' collinearity with their destination:
	
	 {\bf Case 1-Collinear Agents:} Suppose that agents $i$ and $j$ are collinear with the destination at the contact instance $t=t^f_{ij}$, and  that agent $i$ stays in front of agent $j$  relative to the destination. Then $\al_{i}(t^f_{ij}) = \al_{j}(t^f_{ij}) = \al_{ij}(t^f_{ij})$. Combing this together with~\eqref{e20} and~\eqref{e21} we come up to 
	\begin{equation}\label{e22}
		\lm\tau s_i\oa_j = \lm \tau s_j\oa_i,
	\end{equation}
	which implies that $s_i\oa_j(t) =  s_j\oa_i(t)$ assuming $\lm>0$; otherwise we do not have enough information to proceed. This scenario, especially prevalent in corridor settings, simplifies the control strategy to maintaining the same velocity for both agents during their contact period; (see the proof in Section~4.1).
	
 {\bf Case 2-Non-Collinear Agents:} In scenarios where agents are not collinear with the destination at contact time, the control strategy becomes more complex. The relationship expressed in equation~\eqref{e20} does not provide sufficient information to deduce a proportional relationship as in equation~\eqref{e22}. Future work may focus on exploring more intricate control laws that can handle the nuanced interactions. 

\subsection{The Controlled Crowd Motion Models with Multiple Agents in a Corridor}
In this subsection, we examine the controlled crowd motion model for multiple agents in a corridor setting, employing specific assumptions to simplify the complex interactions. We assume that all agents are oriented such that the destination is always directly to their right. This assumption leads to the agent having fixed angles 
	$$\begin{cases}
		\al_i(t) = 180^\circ, \mbox{ for all } i=1, 2, \ldots, n, \\
		\al_{i i+1}(t) = 0^\circ, \mbox{for all } i =1, 2, \ldots, n-1,
	\end{cases}	
	$$
	significantly simplifying the analysis. 
	Then the differential relation in~\eqref{e18a} can be read as 
	\begin{equation}
		\label{e23}
		\begin{cases}
			\dox_1(t) = (\oa_1s_1, 0) - (\eta_{12}(t), 0), \\
 			\vdots\\
			\dox_j(t) =  (\oa_js_j, 0) + (\eta_{j-1j}(t)-\eta_{j+1j}(t), 0), \mbox{ for } j=2, \ldots, n-1,\\
			\vdots \\
			\dox_n(t) = (\oa_n s_n, 0) + (\eta_{n-1n}(t), 0),
		\end{cases}
	\end{equation} 
	which provides a clear picture of the agent's velocities, taking into account the interactions among adjacent agents. For example, agent $j$ adjusts its velocity based on the interactions with agents $j-1$ and $j+1$.

 The relationship between the control functions $\oa_i(\cdot)$ and the vectors $\q^x(\cdot)$ in~\eqref{e21} is reflected by $\lm\tau \oa_i = -s_iq^x_{i1}(t)$
for a.e. $t\in [0,T]$ and $i=1,2,\ldots, n$, leading us to understand how agents in contact synchronize their velocities. Specifically, if agents $i$ and $i+1$ are in contact, meaning that $\eta_{i,i+1}(t)>0$, equations~\eqref{e20} and~\eqref{e22} can be respectively read as 
$q^x_{i1}(t) = q^x_{i+1,1}(t)$
and 
$\lm\tau s_i\oa_{i+1} = \lm\tau s_{i+1}\oa_i $,
which implies that
$s_i\oa_{i+1} = s_{i+1}\oa_i$
under the assumption that $\lm>0$. 

 During contact times $[t^f_{ii+1}, t^l_{ii+1}]$, agents $i$ and $i+1$ maintain a constant distance, as indicated as follows
\begin{equation}
	\label{rel1}
	\| \bar\x_{i+1}-\bar\x_i(t) \|^2 =(L_i+l_{i+1})^2 \mbox{ for all } t\in [t^f_{ii+1}, t^l_{ii+1}].
\end{equation}
Differentiating equation~\eqref{rel1} with respect to $t$ gives 
$$
\la \bar\x_{i+1}(t)-\bar\x_i(t), \dox_{i+1}(t)-\dox_i(t)\ra=0, 
$$
or, equivalently 
\begin{equation}
	\label{rel2}
[\ox_{i+1,1}(t)-\ox_{i,1}(t)][\dot\ox_{i+1,1}(t)-\dot\ox_{i,1}(t)] = 0,
\end{equation}
since the second components of vectors $\bar\x_i(t)$ and $\dox_i(t)$ are identical to zero in the corridor settings. It then follows from~\eqref{rel2} that $\dot\ox_{i+1,1}(t)-\dot\ox_{i,1}(t)=0$ and thus $\dox_{i+1}(t)-\dox_{i}(t)$ for all $t\in [t^f_{ii+1}, t^l_{ii+1}]$, and hence the velocities of the agents must be identical during their period of contact. 

 In conclusion, these simplifications and deductions allow us to understand how multiple agents in a corridor setting adjust their velocities in response to each other. When in contact, adjacent agents synchronize their speeds to ensure a consistent and coordinated movement towards the destination. This insight is crucial for developing effective control strategies in such constrained environments, where agents need to navigate efficiently while maintaining a safe distance from each other.
 
\subsection{The Crowd Motion Model with Two Agents}
	In this section, we explore a model involving two agents in a corridor, focusing on their velocity dynamics. Building on equation~\eqref{e23} we establish a relationship between their velocities as	
	\begin{equation}
		\label{e24}
		\begin{cases}
			\dox_1(t) = (\oa_1s_1 - \eta_{12}(t), 0),\\
			\dox_2(t) = (\oa_2s_2 + \eta_{12}(t), 0).
		\end{cases}
	\end{equation}
	We consider different scenarios where agents interact and adjust their velocities, categorized into distinct cases.
	
	 {\bf Case 1: }{\it Agents in contact} \\
	In this scenario, the agents are in direct contact for the time interval $[t^f_{12}, T]$. They synchronize their velocities $\dox_1(t)=\dox_2(t)$ during this period. This synchronization allows us to determine the interaction force $\eta_{12}(t)$ as follows
	\begin{equation}
		\label{e25}
		\begin{array}{ll}
		\eta_{12}(t) = \frac{1}{2}(\oa_1s_1-\oa_2s_2)>0
		\end{array}
	\end{equation}
	for all $t\in[t^f_{12},T]$. Moreover, considering the interaction equation from ~\eqref{e22}, we can explicitly compute $\eta_{12}(t)$ as 
	\begin{equation}
		\label{e26}
		\eta_{12}(t) = 
		\begin{cases}
			0& \mbox{ if } t\in [0,t^f_{12}),\\
			\frac{\oa_2(s_1^2-s^2_2)}{2s_2}& \mbox{ if } t\in [t^f_{12},T],
		\end{cases}
	\end{equation}
	thanks to~\eqref{e25}. This formulation implies $s_1>s_2$ due to the relative distances of agents from the destination. 
The interaction between the agents adjusts their speed to maintain contact, with agent 1 decelerating and agent 2 accelerating by the quantity $\eta_{12}(t)$.

\subsubsection*{Trajectory Calculations}\vspace{-0.5cm}
Using the Newton-Leibniz formula, the trajectories of the agents are derived from~\eqref{e24}, considering their synchronized movement since the time $t=t^f_{12}$: 
	\begin{equation}
		\label{e26a}
		\begin{cases}
			\bar\x_1(t) = \bar\x_1(0) + (t\oa_1s_1-t\eta_{12}(t^f_{12})+t^f_{12}\eta_{12}(t^f_{12}),0), \\
			\bar\x_2(t) = \bar\x_2(0) + (t\oa_2s_2+t\eta_{12}(t^f_{12})-t^f_{12}\eta_{12}(t^f_{12}),0),
		\end{cases}
	\end{equation}
	for all $t\in [t^f_{12}, T]$. 
	
	\subsubsection*{Contact Time and Interaction Effort}\vspace{-0.5cm}
	The contact time $t^f_{12}$ and the interaction effort are deduced from the agent's initial positions and the constraint that their separation remains constant, i.e. $\l\|\bar\x_2(t) - \bar\x_1(t) \r\| = L_2+L_1$:
	\begin{equation}
		\label{e27} 
		\begin{cases}
			t^f_{12}\eta_{12}(t^f_{12}) = \frac{1}{2}[\ox_{21}(0) - \ox_{11}(0) - (L_1 + L_2)] = \Lambda_{12},\\
			t^f_{12} = \frac{\Lm_{12}}{\eta_{12}(t^f_{12})}. 
		\end{cases}
	\end{equation}

\subsubsection*{Cost Functional}\vspace{-0.5cm}
	The cost functional, influenced by the control variable $\oa_2$, can be computed as
	\begin{equation}
		\label{e28}
		\begin{array}{ll}
		J[\bar\x,\oa] &= \frac{1}{2}\l[ \l\| \bar\x_1(T) - \xd\r\|^2 +  \l\| \bar\x_2(T) - \xd\r\|^2 \r] + \frac{\tau T}{2}(\oa_1^2+\oa_2^2).
		\end{array}
	\end{equation}
	We further examine cases where agents start in contact (Case 1a) or apart (Case 1b), each affecting their movement and the cost functional. These scenarios are critical for understanding how initial conditions influence their paths and interaction dynamics.
	
\noindent	 {\bf Case 1a:} {\it The agents are initially in contact, i.e. $\bar\x_{21}(0)-\bar\x_{11}(0)=L_1+L_2$}. In this very situation, we have $t^f_{12}=0$ and the agents are heading to the destination with the same velocity until the end of the process. The cost functional is computed as 
	\begin{equation}
		\label{e30}
		\begin{array}{ll}
		J[\bar\x,\oa] &=  \frac{1}{2}\bigg[ \l\vert \bar\x_{11}(0) +T\frac{\oa_2(s_1^2+s^2_2)}{2s_2} - \xd_1\r\vert^2 \\[1ex]
		&+  \l\vert \bar\x_{21}(0) + T\frac{\oa_2(s_1^2+s^2_2)}{2s_2}- \xd_1\r\vert^2 \bigg] + \frac{\tau T\oa_2^2(s_1^2+s^2_2)}{2s_2^2}.
		\end{array}
	\end{equation}
	Minimizing this quadratic function with respect to $\oa_2$ gives an optimal pair $(\oa_1,\oa_2)$ to our optimal control problem. 
	
	\noindent{\bf Case 1b:} {\it The agents are out of contact at the beginning i.e. $t^f_{12}>0.$} It then follows from~\eqref{e26a},~\eqref{e27} and the fact $t^f_{12}\leq T$ that 
	\begin{equation}
		\label{e31}
		\begin{array}{ll}
		\frac{\Lm_{12}}{\eta_{12}(t^f_{12})}\leq T \iff \eta_{12}(t^f_{12})\geq \frac{\Lm_{12}}{T} \iff \oa_2\geq \frac{2s_2\Lm_{12}}{T(s^2_1-s^2_2)}
		\end{array}
	\end{equation}
	The cost functional in this case is 
	\begin{equation}
		\label{e32}
		\begin{array}{ll}
			J[\bar\x,\oa]  &=\frac{1}{2}\bigg[ \l\| \bar\x_{11}(0) +T\oa_1s_1-T\eta_{12}(t^f_{12})+t^f_{12}\eta_{12}(t^f_{12}) - \xd_1\r\|^2\\
			&+  \l\| \bar\x_{21}(0) + T\oa_2s_2+T\eta_{12}(t^f_{12})-t^f_{12}\eta_{12}(t^f_{12})- \xd_1\r\|^2 \bigg] \\
			&+ \frac{\tau T\oa_2^2(s_1^2+s^2_2)}{2s_2^2}.
		\end{array}
	\end{equation}
	Minimizing this quadratic function with respect to $\oa_2$ subject to constraint~\eqref{e31} gives an optimal pair $(\oa_1,\oa_2)$ to our optimal control problem. 
	
	\noindent{\bf Case 2:} {\it No Contact Between Agents Throughout}\\
	In this scenario, the agents do not come into contact at any point during the interval $[0,T]$. Consequently, the interaction term $\eta_{12}(t)$ is inactive, remaining zero throughout this period. This situation simplifies the cost functional significantly:	
	\begin{equation}
		\label{e33}
		\begin{array}{ll}
		J[\bar\x,\oa] &=  \frac{1}{2}\l[ \l\vert \bar\x_{11}(0) + T\oa_1s_1 - \xd_1\r\vert^2 +  \l\vert \bar\x_{21}(0) + T\oa_2s_2 - \xd_1\r\vert^2 \r] \\
		&+ \frac{\tau T}{2}(\oa_1^2+\oa_2^2).
		\end{array}
	\end{equation}
	This functional represents the scenario where each agent independently navigates towards the destination, unaffected by the other's presence. The distance between them at any time $y$ is given by
	$$
	\l\| \x_2(t) - \x_1(t)\r\| = \l| \ox_{21}(0) - \ox_{11}(0) + t(\oa_2s_2 - \oa_1s_1)\r|> L_1+L_2,
	$$
	ensuring no contact throughout their journey. The inequality condition for this non-contact scenario is represented as:
	\begin{equation}
		\label{e34}
		t(\oa_1s_1 - \oa_2s_2) < \ox_{21}(0) - \ox_{11}(0) - (L_1+L_2)
	\end{equation}
	for all $t\in [0,T]$.
	
	 $\bullet$ If $\oa_1s_1 < \oa_2s_2$: this case naturally satisfies the non-contact condition, implying agent 2 moves faster than agent 1, ensuring they remain apart.
	
	 $\bullet$ If $\oa_1s_1 > \oa_2s_2$: here, the inequality in equation~\eqref{e34} sets an upper limit for the time period during which the agents can stay out of contact, based on their initial positions and velocities. 

 The work presented in the previous sections is encapsulated in  Algorithm~1. This algorithm can be seen as an extension of the methodologies discussed in~\cite{cm4}, adapted to more general data settings. To demonstrate the effectiveness of this algorithm, we present two examples with distinct datasets. 

\begin{figure}[!ht]
	\centering
	\includegraphics[scale=0.4]{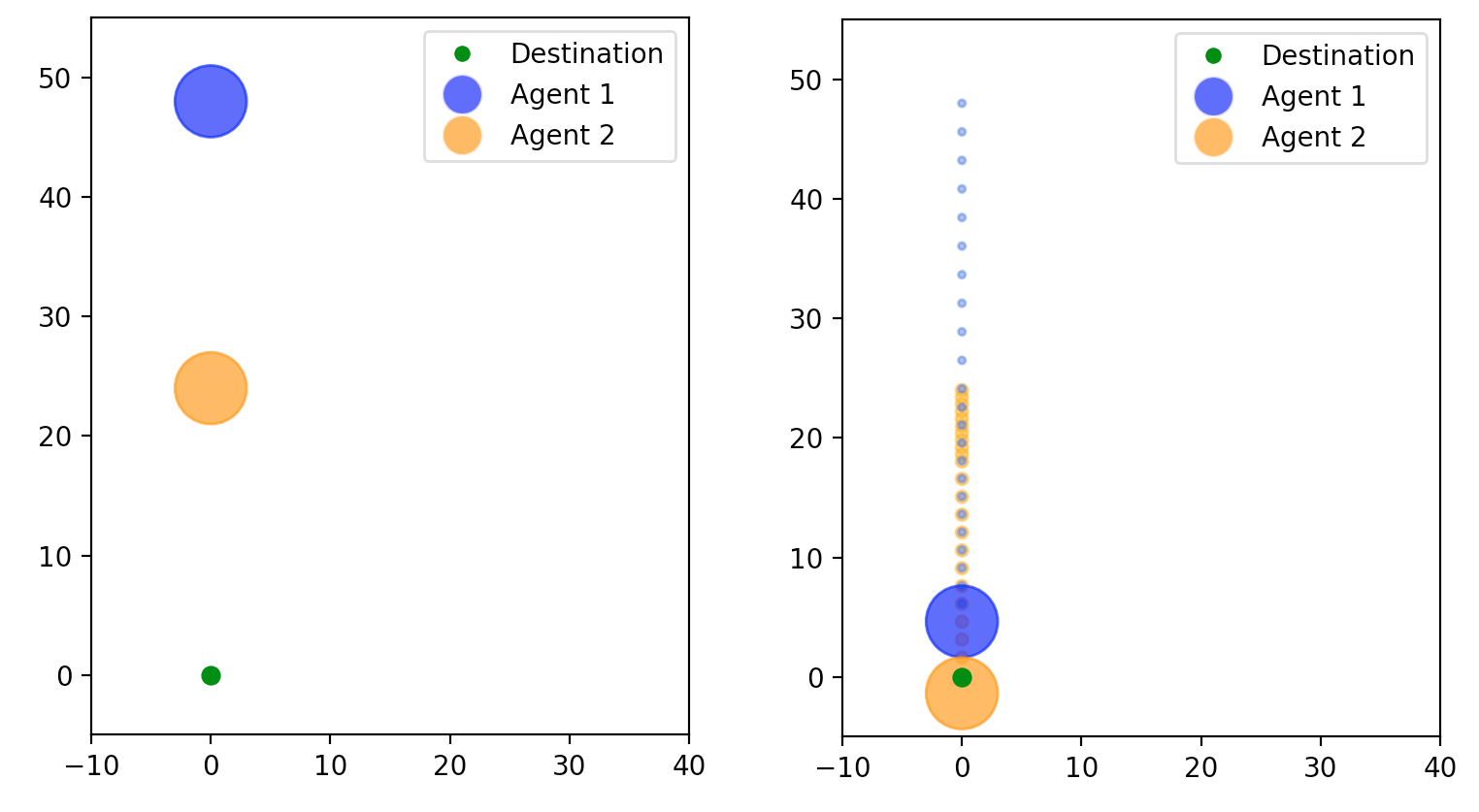}  
	\caption{Illustration for Example~\ref{eg51} at $t=0$ and $t=T$}\label{eg51fig}
\end{figure}

\begin{example}\label{eg51} {\bf Basic Crowd Motion Scenario} \\
Consider a controlled crowd motion problem characterized by
$$
\begin{cases}
	\xd=(0,0), \; \x_1^0 = (0,48), \; \x_2^0= (0,24), \\
	T = 6, \; L_1=L_2=3.
\end{cases}
$$
This example demonstrates the agents' performance under varying values of the parameter $\tau$. The results are tabulated below and further illustrated in Figure~\ref{eg51fig}, showing agent positions at $t=0$ and $t=T$ . 
	{\small \begin{center}
			\begin{tabular}{|c|c|c|c|c|}
				\hline
				$\tau$ &         $\oa_1$ &        $\oa_2$ &        $t^f_{12}$ &        $J[\bar{x},\bar{a}]$ \\
				\hline
				1.0 &  1.195021 &  0.597510 &  2.510417 &  14.377593 \\
				2.0 &  1.190083 &  0.595041 &  2.520833 &  19.710744 \\
				3.0 &  1.185185 &  0.592593 &  2.531250 &  25.000000 \\
				4.0 &  1.180328 &  0.590164 &  2.541667 &  30.245902 \\
				5.0 &  1.175510 &  0.587755 &  2.552083 &  35.448980 \\
				6.0 &  1.170732 &  0.585366 &  2.562500 &  40.609756 \\
				7.0 &  1.165992 &  0.582996 &  2.572917 &  45.728745 \\
				8.0 &  1.161290 &  0.580645 &  2.583333 &  50.806452 \\
				9.0 &  1.156626 &  0.578313 &  2.593750 &  55.843373 \\
				10.0 &  1.152000 &  0.576000 &  2.604167 &  60.840000 \\
				\hline
			\end{tabular}
	\end{center}}
\end{example}

\begin{figure}[!ht]
	\centering
	\includegraphics[scale=0.3]{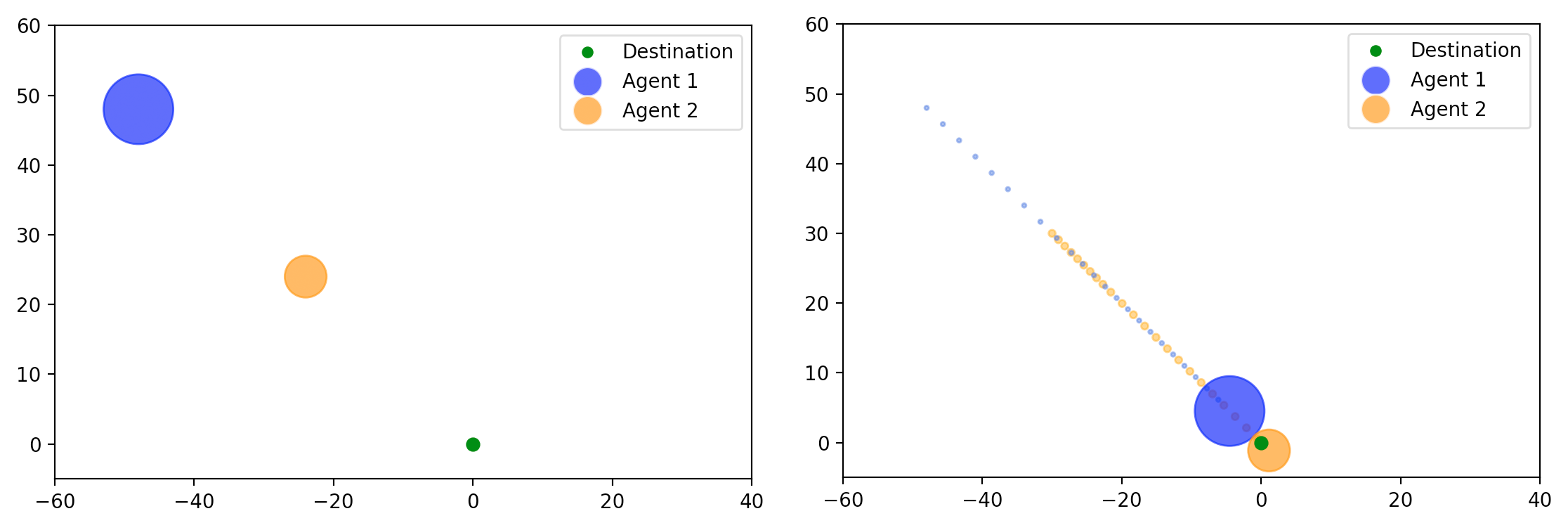}  
	\caption{Illustration for Example~\ref{eg52} at $t=0$ and $t=T$}\label{eg52fig}
\end{figure}

\begin{example}\label{eg52}{\bf Advanced Crowd Motion Scenario}\\
Here, we explore a more complex crowd motion setting:
$$
\begin{cases}
	\xd=(0,0), \; \x_1^0 = (-48,48), \; \x_2^0= (-24,24), \\
	T = 6, \; L_1=5, \; L_2=3.
\end{cases}
$$
The optimal values of $\oa$ and $J[\bar\x,\oa]$ are computed for different $\tau$ values, showing improved agent performance compared to uncontrolled scenarios. The results are summarized in the table below, and the agents' trajectories are depicted in Figure~\ref{eg52fig}.
	{\small\begin{center}
		\begin{tabular}{|c|c|c|c|c|}
			\hline
			$\tau$ &         $\oa_1$ &        $\oa_2$ &        $t^f_{12}$ &           $J[\bar{x},\bar{a}]$ \\
			\hline
			1.0 &  1.166355 &  0.728972 &  2.170803 &  21.685981 \\
			2.0 &  1.164179 &  0.727612 &  2.174861 &  27.350746 \\
			3.0 &  1.162011 &  0.726257 &  2.178918 &  32.994413 \\
			4.0 &  1.159851 &  0.724907 &  2.182976 &  38.617100 \\
			5.0 &  1.157699 &  0.723562 &  2.187033 &  44.218924 \\
			6.0 &  1.155556 &  0.722222 &  2.191091 &  49.800000 \\
			7.0 &  1.153420 &  0.720887 &  2.195149 &  55.360444 \\
			8.0 &  1.151291 &  0.719557 &  2.199206 &  60.900369 \\
			9.0 &  1.149171 &  0.718232 &  2.203264 &  66.419890 \\
			10.0 &  1.147059 &  0.716912 &  2.207321 &  71.919118 \\
			\hline
		\end{tabular}
	\end{center}}
The performances of the agents are significantly better than the uncontrolled cases (with $\oa_1=\oa_2=1$), which is clearly shown in the following table.
	{\small\begin{center}
		\begin{tabular}{|c|c|c|}
			\hline
			$\tau$ &    $t^f_{12}$ & $J[\bar{x}]$ \\
			\hline
			1.0 &  4.114382 &  114.16 \\
			2.0 &  4.114382 &  120.16 \\
			3.0 &  4.114382 &  126.16 \\
			4.0 &  4.114382 &  132.16 \\
			5.0 &  4.114382 &  138.16 \\
			6.0 &  4.114382 &  144.16 \\
			7.0 &  4.114382 &  150.16 \\
			8.0 &  4.114382 &  156.16 \\
			9.0 &  4.114382 &  162.16 \\
			10.0 &  4.114382 &  168.16 \\
			\hline
		\end{tabular}
	\end{center}}
	These examples highlight the algorithm's efficacy in optimizing the movement and interaction of agents in varying crowd motion scenarios.
\end{example}

\begin{remark}

The relationship between control parameters $\oa_1$ and $\oa_2$, as delineated in equation~\eqref{e22}, plays a pivotal role in computational efficiency. This interrelation, stemming from the necessary optimality conditions outlined in Theorem~\ref{Th3}, significantly simplifies the computational process. While it is possible to solve the problem by independently minimizing the cost functional with respect to $\oa_1$ and $\oa_2$, this approach drastically increases the complexity of Algorithm~1. Furthermore, this interrelation between controls proves even more beneficial in scenarios involving three or more agents, where computational demands escalate rapidly compared to two-agent cases. To further elucidate this, let us consider the adjoint arcs $\p=(\p^x,\p^a)$ and $\q=(\q^x,\q^a)$ from our necessary optimality conditions in Theorem~\ref{Th3}: 
	\begin{enumerate}
		\item We have $\p^x(t) = \p^x(T)$ and $\p^a(t) ={\bf 0}$ for all $t\in [0,T]$;
		\item The function $\q^x(t)$ is defined as follows:
		$$
		\q^x(t) = 
		\begin{cases}
			\p^x(t) + \gg([t^f_{12},T]) & \mbox{ if } t \in [0, t^f_{12})\\
			\p^x(t) + \gg([t,T]) & \mbox{ if } t \in [t^f_{12},T];
		\end{cases} 
		$$
		\item We observe that $\q^x_1(t) = \q^x_2(t)= \l( -\frac{\lm\tau \oa_1}{s_1},0 \r)= \l(-\frac{\lm\tau \oa_2}{s_2},0\r) $ for a.e. $t\in [0,T]$, that is $\q^x(\cdot)$ takes a constant value on $[0,T]$ and so does the measure $\gg([\cdot,T])$;
		\item The terminal conditions for $\p^x(T)$ are given by
		$$
		\begin{cases}
		\p^x_1(T)= (\lm|\bar\x_{11}(T)-\xd_1| - \eta_{21}(T),0) \\
		\p^x_2(T) = (\lm|\bar\x_{21}(T)-\xd_1| + \eta_{21}(T),0).
		\end{cases}
		$$
	\end{enumerate}
These observations suggest that the measure $\gg([\cdot,T]) $	and the scalar $\lm$ satisfy a system of equations, offering numerous choices for their values.  
 Note that the selection of $\lm , \p, \q,$ and $\gg$ as  $\lm=0, \q(\cdot)\equiv {\bf 0}, \gg = \delta_{\{T\}}(\eta_{21}(T),-\eta_{21}(T)), \p(\cdot) \equiv - (\eta_{21}(T),-\eta_{21}(T)),$ and $\eta_{21}(T)=\frac{\oa_2(s_1^2-s^2_2)}{2s_2}$ leads to a degeneration of the necessary optimality conditions, rendering them less informative for finding optimal solutions. On the other hand, assuming $\lm>0$ allows us to derive a valuable relationship between $\oa_1$ and $\oa_2$ using equation~\eqref{e22}, which is instrumental in computing optimal controls. This underscores the effectiveness of our necessary optimality conditions in addressing complex controlled crowd motion models. 
\end{remark}
	
\subsection{The Crowd Motion Model with Three Agents}\vspace{-0.5cm}

We now extend our analysis to a crowd motion control problem involving three  agents $\bar\x_1, \bar\x_2,$ and $ \bar\x_3$. The differential relationships, adapted from equation~\eqref{e23} are as follows
	\begin{equation}
		\label{e29s}
		\begin{cases}
			\dox_1(t) = (\oa_1s_1 - \eta_{12}(t), 0), \\
			\dox_2(t) =  (\oa_2s_2 + \eta_{12}(t)-\eta_{32}(t), 0),\\
			\dox_3(t) = (\oa_3 s_3 +\eta_{23}(t), 0)
		\end{cases}
	\end{equation} 
	for a.e. $t\in [0,T]$. 
	
\subsubsection*{Agent Interactions}
\vspace{-0.5cm}
	\begin{figure}[!ht]
		\centering
		\includegraphics[scale=.2]{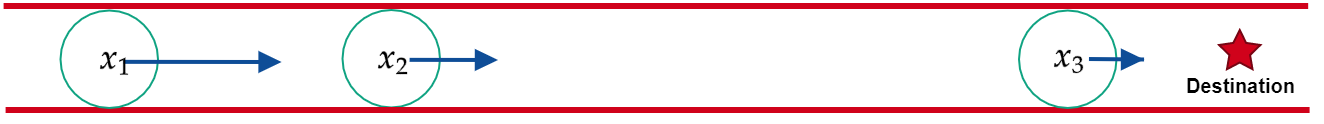}  
		\caption{Three Agents in a Corridor}\label{fig16}
	\end{figure}
	Agent 2, positioned in the middle, must consider both agents 1 and 3. The interaction between any two agents is reciprocal, denoted as $\eta_{32}(\cdot)=\eta_{23}(\cdot)$ in equation~\eqref{e29s}. Denoting $t^f_{123}$ as the first contact time among three agents, we have $\max\l\{t^f_{12}, t^f_{23}\r\}\leq t^f_{123}$.
	
	\subsubsection*{Agent Behaviors at Contact Times}\vspace{-0.5cm}
	From equation~\eqref{e22}, the relationships during contact times are  $\oa_1 = \frac{s_1}{s_2}\oa_2$ and $\oa_2= \frac{s_2}{s_3}\oa_3$ , leading to 
	\begin{equation}
		\label{e30s}
		\begin{array}{ll}
		\oa_2= \frac{s_2}{s_3}\oa_3, \; \oa_1 = \frac{s_1}{s_3}\oa_3
		\end{array}
	\end{equation}
	at $t=t^f_{123}$. These relationships are crucial for understanding the agent's velocity adjustments to maintain contact. \vspace{-0.5cm}
	
	\subsubsection*{Interaction Efforts at Contact Times}\vspace{-0.5cm}
	Considering the velocities at various contact times, we derive the following relationships for interaction efforts:
	\begin{equation}
		\label{e31s}
		\begin{cases}
			\eta_{12}(t^f_{12}) = \frac{1}{2}\l[\oa_1s_1-\oa_2s_2 + \eta_{23}(t^f_{12}) \r]\\[1ex]
			\eta_{23}(t^f_{23}) =  \frac{1}{2}\l[\oa_2s_2-\oa_3s_3 + \eta_{12}(t^f_{23}) \r].
		\end{cases}
	\end{equation}
Furthermore, the equality of velocities when all three agents are in contact $\dox_1(t^f_{123}) = \dox_2(t^f_{123}) = \dox_3(t^f_{123})$ yields
	$$
	\begin{cases}
		\eta_{12}(t^f_{123}) = \frac{1}{3}(2\oa_1s_1-\oa_2s_2-\oa_3s_3) \\
		\eta_{23}(t^f_{123}) = \frac{1}{3}(\oa_1s_1 +\oa_2s_2-2\oa_3s_3).
	\end{cases}
	$$
	Combining these relations with equations in~\eqref{e30s} we can express $\eta_{12}(t^f_{123})$ and $\eta_{23}(t^f_{123})$ in terms of $\oa_3$ as follows 
	\begin{equation*}
		\label{e32s}
		\begin{cases}
			\eta_{12}(t^f_{123}) = \frac{(2s^2_1-s^2_2-s^2_3)\oa_3}{3s_3} \\
			\eta_{23}(t^f_{123}) = \frac{(s^2_1+s^2_2-2s^2_3)\oa_3}{3s_3}.
		\end{cases}
	\end{equation*}
	
	\subsubsection*{Velocity and Trajectory Calculations}\vspace{-0.5cm}
The velocities of three agents during the time interval $[t^f_{123}, T]$ are thus	
	$$
	\dox_1(t) = \dox_2(t) = \dox_3(t) = \l(\frac{(s^2_1+s^2_2+s^2_3)\oa_3}{3s_3},0\r),\quad \mbox{ for all } t\in [t^f_{123}, T].
	$$
The trajectories and velocities are then determined for various cases, including when agents 1 and 2 contact first (Case 1), when agents 2 and 3 contact first, and when all three agents contact simultaneously (Case 3).
	
\noindent$\bullet$  {\bf Case 1: Agent 1 and 2 contact first}, i.e. $t^f_{12} < t^f_{23},$.  In this case there is no contact between agent 2 and 3 at $t^f_{12}$, i.e. $\eta_{23}(t^f_{12}) = 0$, and so $t^f_{23}=t^f_{123}$. Thus, the first equation in~\eqref{e31s} implies 
	$
	\eta_{12}(t^f_{12}) = \frac{1}{2}(\oa_1s_1-\oa_2s_2)>0,
	$
	which reduces to~\eqref{e25}. The scalar functions $\eta_{12}(\cdot)$ and $\eta_{23}(\cdot)$ showing the interaction efforts among three agents are given by
	\begin{equation}
		\label{e34s}
		\eta_{12}(t) = 
		\begin{cases}
			0 &\mbox{ if }	t \in [0,t^f_{12}) \\
			\frac{1}{2}(\oa_1s_1-\oa_2s_2) = \frac{(s_1^2-s^2_2)\oa_2}{2s_2} &  \mbox{ if } t \in [t^f_{12}, t^f_{123}) \\
			\frac{1}{3}(2\oa_1s_1-\oa_2s_2-\oa_3s_3)= \frac{(2s^2_1-s^2_2-s^2_3)\oa_3}{3s_3} & \mbox{ if } t \in [t^f_{123},T]
		\end{cases}
	\end{equation}
	and 
	\begin{equation}
		\label{e35s}
		\eta_{23}(t) = 
		\begin{cases}
			0 &\mbox{ if }	t \in [0,t^f_{12}) \\   
			0 &  \mbox{ if } t \in [t^f_{12}, t^f_{123}) \\
			\frac{1}{3}(\oa_1s_1 +\oa_2s_2-2\oa_3s_3) =  \frac{(s^2_1+s^2_2-2s^2_3)\oa_3}{3s_3} & \mbox{ if } t \in [t^f_{123},T].
		\end{cases}
	\end{equation}
Let us next compute  the contact time $t^f_{12}$ for agents 1 and 2 and $t^f_{23}$ for agents 2 and 3 respectively. Indeed, thanks to~\eqref{e34s}, we have
	$$
	\begin{array}{ll}
	\bar\x_2(t^f_{12}) - \bar\x_1(t^f_{12})&= \bar\x_2(0) - \bar\x_1(0) + \l(t^f_{12}(\oa_2s_2-\oa_1s_1),0\r)\\
  &=\bar\x_2(0) - \bar\x_1(0) + \l(-2t^f_{12}\eta_{12}(t^f_{12}),0 \r) 
 \end{array}
 	$$
 	and
	$$
	\begin{array}{ll}
	\bar\x_3(t^f_{23}) - \bar\x_2(t^f_{23})&= \bar\x_3(0) - \bar\x_2(0) \\&\quad+ \Big(t^f_{23}(\oa_3s_3-\oa_2s_2-\eta_{12}(t^f_{12})) + t^f_{12}\eta_{12}(t^f_{12}), 0 \Big)\\
		&= \bar\x_3(0) - \bar\x_2(0) + \l( t^f_{23}\frac{2\oa_3s_3-\oa_2s_2-\oa_1s_1}{2}+t^f_{12}\eta_{12}(t^f_{12}),0\r).
	\end{array}
	$$
It then follows from $\l\| \bar\x_2(t^f_{12}) - \bar\x_1(t^f_{12})\r\| = L_1+L_2$ and $\l\| \bar \x_3(t^f_{23}) - \bar\x_1(t^f_{23})\r\| = L_2+L_3$ that 
	\begin{equation}
		\label{e31a}
		\begin{cases}
			t^f_{12}\eta_{12}(t^f_{12}) &= \frac{\ox_{21}(0)-\ox_{11}(0)-(L_1+L_2)}{2}: = \Lm_{12},\\
			t^f_{12}& = \frac{\Lm_{12}}{\eta_{12}(t^f_{12})}, \\  
			\frac{3}{2}t^f_{23}\eta_{23}(t^f_{23}) &=  2\Lm_{23}+\Lm_{12},\\ 
			t^f_{123} = t^f_{23} & =  \frac{2(2\Lm_{23}+\Lm_{12})}{3\eta_{23}(t^f_{123})}\leq T,
		\end{cases}
	\end{equation}
where $\Lm_{23} := \frac{\ox_{31}(0)-\ox_{21}(0)-(L_2+L_3)}{2}$. Combing~\eqref{e34s},~\eqref{e35s}, and~\eqref{e31a} enables us to express the contact time $t^f_{12}$ and $t^f_{23} = t^f_{123}$ in terms of the controls $\oa_1, \oa_2$ and $\oa_3$. Hence, the velocities and the corresponding trajectories of three agents can be computed in terms of the controls respectively as follows 
$$
		\dox_1(t) = 
		\begin{cases}
			(\oa_1s_1, 0) & \mbox{ if } t \in [0,t^f_{12}) \\
			\l(\frac{(s^2_1+s^2_2)\oa_2}{2s_2},0 \r)  &  \mbox{ if } t \in [t^f_{12}, t^f_{123}) \\
			\l(\frac{(s^2_1+s^2_2+s^2_3)\oa_3}{3s_3},0 \r) & \mbox{ if } t \in [t^f_{123},T],
		\end{cases}
$$
$$
		\dox_2(t) = 
		\begin{cases}
			(\oa_2s_2, 0) & \mbox{ if } t \in [0,t^f_{12}) \\
			\l(\frac{(s^2_1+s^2_2)\oa_2}{2s_2},0 \r)  &  \mbox{ if } t \in [t^f_{12}, t^f_{123}) \\
			\l(\frac{(s^2_1+s^2_2+s^2_3)\oa_3}{3s_3},0 \r) & \mbox{ if } t \in [t^f_{123},T],
		\end{cases}
$$
$$
		\dox_3(t) = 
		\begin{cases}
			(\oa_3s_3, 0) & \mbox{ if } t \in [0,t^f_{12}) \\
			(\oa_3s_3, 0)  &  \mbox{ if } t \in [t^f_{12}, t^f_{123}) \\
			\l(\frac{(s^2_1+s^2_2+s^2_3)\oa_3}{3s_3},0 \r) & \mbox{ if } t \in [t^f_{123},T],
		\end{cases}
$$
	and 
	\begin{equation}
		\label{e37a}
		\bar\x_1(t) = 
		\begin{cases}
			\bar \x_1(0) + (t\oa_1s_1, 0) & \mbox{ if } t \in [0,t^f_{12}) \\
			\bar\x_1(0) + \l( t^f_{12}\oa_1s_1 + (t-t^f_{12})\frac{(s^2_1+s^2_2)\oa_2}{2s_2},0 \r)  &  \mbox{ if } t \in [t^f_{12}, t^f_{123}) \\
			\bar\x_1(0) + \bigg(  t^f_{12}\oa_1s_1 + (t^f_{123}-t^f_{12})\frac{(s^2_1+s^2_2)\oa_2}{2s_2} & \\
			 \quad + (t-t^f_{123})\frac{(s^2_1+s^2_2+s^2_3)\oa_3}{3s_3},0 \bigg) & \mbox{ if } t \in [t^f_{123},T],
		\end{cases}
	\end{equation}
	\begin{equation}
		\label{e37b}
		\bar\x_2(t) = 
		\begin{cases}
			\bar \x_2(0) + (t\oa_2s_2, 0) & \mbox{ if } t \in [0,t^f_{12}) \\
			\bar\x_2(0) + \l( t^f_{12}\oa_2s_2 + (t-t^f_{12})\frac{(s^2_1+s^2_2)\oa_2}{2s_2},0 \r)  &  \mbox{ if } t \in [t^f_{12}, t^f_{123}) \\
			\bar\x_2(0) + \bigg(  t^f_{12}\oa_2s_2 + (t^f_{123}-t^f_{12})\frac{(s^2_1+s^2_2)\oa_2}{2s_2}\\
			 +(t-t^f_{123})\frac{(s^2_1+s^2_2+s^2_3)\oa_3}{3s_3},0 \bigg) & \mbox{ if } t \in [t^f_{123},T],
		\end{cases}
	\end{equation}
	\begin{equation}
		\label{e37c}
		\bar\x_3(t) = 
		\begin{cases}
			\bar \x_3(0) + (t\oa_3s_3, 0) & \mbox{ if } t \in [0,t^f_{12}) \\
			\bar\x_3(0) +(t\oa_3s_3, 0)  &  \mbox{ if } t \in [t^f_{12}, t^f_{123}) \\
			\bar\x_3(0) + \l(  t^f_{123}\oa_3s_3  +(t-t^f_{123})\frac{(s^2_1+s^2_2+s^2_3)\oa_3}{3s_3},0 \r) & \mbox{ if } t \in [t^f_{123},T].
		\end{cases}
	\end{equation}
	Thus the cost functional given by 
	\begin{equation}
		\label{e38}
		J[\bar\x,\oa] = \frac{1}{2}\sum^3_{i=1} \l\| \bar\x_i(T) - \xd\r\|^2  + \dfrac{\tau T}{2}(\oa_1^2+\oa_2^2+\oa_3^2)
	\end{equation}
	can be expressed in terms of $\oa_3$. 
	
	\noindent$\bullet$ {\bf Case 2: Agent 2 and 3 contact first}, i.e. $t^f_{23} < t^f_{12}$. Using the similar arguments as in case 1 allows us to express the interaction efforts among three agents, their contact times, velocities, trajectories, and the cost functional in terms of the control $\oa_3$.
	The interaction-effort functions $\eta_{12}(\cdot)$ and $\eta_{23}(\cdot)$ are given by 
	\begin{equation}
		\label{e39}
		\eta_{12}(t) = 
		\begin{cases}
			0 &\mbox{ if }	t \in [0,t^f_{23 ,}) \\
			0 &  \mbox{ if } t \in [t^f_{23}, t^f_{123}) \\
			\frac{1}{3}(2\oa_1s_1-\oa_2s_2-\oa_3s_3) =\frac{(2s^2_1-s^2_2-s^2_3)\oa_3}{3s_3} & \mbox{ if } t \in [t^f_{123},T]
		\end{cases}
	\end{equation}
	and 
	\begin{equation}
		\label{e40}
		\eta_{23}(t) = 
		\begin{cases}
			0 &\mbox{ if }	t \in [0,t^f_{23}) \\
			\frac{1}{2}(\oa_2s_2-\oa_3s_3) = \frac{(s^2_2-s^2_3)\oa_2}{2s_2} &  \mbox{ if } t \in [t^f_{23}, t^f_{123}) \\
			\frac{1}{3}(\oa_1s_1 +\oa_2s_2-2\oa_3s_3) = \frac{(s^2_1+s^2_2-2s^2_3)\oa_3}{3s_3} & \mbox{ if } t \in [t^f_{123},T].
		\end{cases}
	\end{equation}
The contact time $t^f_{12}$ and $t^f_{23}$ are as follows
	\begin{equation}
		\label{e31b}
		\begin{cases}
			t^f_{23}\eta_{23}(t^f_{23}) &= \frac{\ox_{31}-\ox_{21}(0)-(L_2+L_3)}{2} = \Lm_{23}\\
			t^f_{23}& = \frac{\Lm_{23}}{\eta_{23}(t^f_{23})} \\
			\frac{3}{2}t^f_{12}\eta_{12}(t^f_{12}) &=  2\Lm_{12}+\Lm_{23}\\
			t^f_{123} = t^f_{12} & =  \frac{2(2\Lm_{12}+\Lm_{23})}{3\eta_{12}(t^f_{12})} \leq T.
		\end{cases}
	\end{equation}
	The velocities and trajectories of the agents are as follows 
$$
		\dox_1(t) = 
		\begin{cases}
			(\oa_1s_1, 0) & \mbox{ if } t \in [0,t^f_{23}) \\
			(\oa_1s_1, 0)  &  \mbox{ if } t \in [t^f_{23}, t^f_{123}) \\
			\l(\frac{(s^2_1+s^2_2+s^2_3)\oa_3}{3s_3},0 \r) & \mbox{ if } t \in [t^f_{123},T],
		\end{cases}
$$
$$
		\dox_2(t) = 
		\begin{cases}
			(\oa_2s_2, 0) & \mbox{ if } t \in [0,t^f_{23}) \\
			\l(\frac{(s^2_2+s^2_3)\oa_3}{2s_3},0 \r)  &  \mbox{ if } t \in [t^f_{23}, t^f_{123}) \\
			\l(\frac{(s^2_1+s^2_2+s^2_3)\oa_3}{3s_3},0 \r) & \mbox{ if } t \in [t^f_{123},T],
		\end{cases}
$$
$$
		\dox_3(t) = 
		\begin{cases}
			(\oa_3s_3, 0) & \mbox{ if } t \in [0,t^f_{23}) \\
			\l(\frac{(s^2_2+s^2_3)\oa_3}{2s_3},0 \r)  &  \mbox{ if } t \in [t^f_{23}, t^f_{123}) \\
			\l(\frac{(s^2_1+s^2_2+s^2_3)\oa_3}{3s_3},0 \r) & \mbox{ if } t \in [t^f_{123},T],
		\end{cases}
$$
	and 
	\begin{equation}
		\label{e42a}
		\bar\x_1(t) = 
		\begin{cases}
			\bar \x_1(0) + (t\oa_1s_1, 0) & \mbox{ if } t \in [0,t^f_{23}) \\
			\bar\x_1(0) + (t\oa_1s_1, 0) &  \mbox{ if } t \in [t^f_{23}, t^f_{123}) \\
			\bar\x_1(0) + \l(  t^f_{123}\oa_1s_1  +(t-t^f_{123})\frac{(s^2_1+s^2_2+s^2_3)\oa_3}{3s_3},0 \r) & \mbox{ if } t \in [t^f_{123},T],
		\end{cases}
	\end{equation}
	\begin{equation}
		\label{e42b}
		\bar\x_2(t) = 
		\begin{cases}
			\bar \x_2(0) + (t\oa_2s_2, 0) & \mbox{ if } t \in [0,t^f_{23}) \\
			\bar\x_2(0) + \l( t^f_{23}\oa_2s_2 + (t-t^f_{23})\frac{(s^2_2+s^2_3)\oa_3}{2s_3},0 \r)  &  \mbox{ if } t \in [t^f_{23}, t^f_{123}) \\
			\bar\x_2(0) + \bigg(  t^f_{23}\oa_2s_2 + (t^f_{123}-t^f_{23})\frac{(s^2_2+s^2_3)\oa_3}{2s_3}&\\
			 +(t-t^f_{123})\frac{(s^2_1+s^2_2+s^2_3)\oa_3}{3s_3},0 \bigg) & \mbox{ if } t \in [t^f_{123},T],
		\end{cases}
	\end{equation}
	\begin{equation}
		\label{e42c}
		\bar\x_3(t) = 
		\begin{cases}
			\bar \x_3(0) + (t\oa_3s_3, 0) & \mbox{ if } t \in [0,t^f_{23}) \\ 
			\bar\x_3(0) + \l( t^f_{23}\oa_3s_3 + (t-t^f_{23})\frac{(s^2_2+s^2_3)\oa_3}{2s_3},0 \r)  &  \mbox{ if } t \in [t^f_{23}, t^f_{123}) \\ 
			\bar\x_2(0) + \bigg(  t^f_{23}\oa_3s_3 + (t^f_{123}-t^f_{23})\frac{(s^2_2+s^2_3)\oa_3}{2s_3}&\\
			 \qquad+(t-t^f_{123})\frac{(s^2_1+s^2_2+s^2_3)\oa_3}{3s_3},0 \bigg) & \mbox{ if } t \in [t^f_{123},T],
		\end{cases}
	\end{equation}
	
	\noindent$\bullet$  {\bf Case 3: All three agents contact simultaneously}, i.e. $t^f_{12} = t^f_{23} = t^f_{123}$. In this scenario, the interaction-effort functions $\eta_{12}(\cdot)$ and $\eta_{23}(\cdot)$ are simply given by 
	\begin{equation}
		\label{e43}
		\eta_{12}(t) = \begin{cases}
			0& \mbox{ if } t\in [0,t^f_{123}), \\
			\frac{(2s^2_1-s^2_2-s^2_3)\oa_3}{3s_3} & \mbox{ if } t \in [t^f_{123}, T]
		\end{cases}
	\end{equation}
	and 
	\begin{equation}
		\label{e44}
		\eta_{23}(t) = \begin{cases}
			0& \mbox{ if } t\in [0,t^f_{123}), \\
			\frac{(s^2_1+s^2_2-2s^2_3)\oa_3}{3s_3} & \mbox{ if } t \in [t^f_{123}, T]. 
		\end{cases}
	\end{equation}
	The agents' velocities and trajectories of the agents are as follows 
	\begin{equation*}
		\label{e45}
		\dox_i(t) = 
		\begin{cases}
			(\oa_is_i, 0) & \mbox{ if } t \in [0,t^f_{123}) \\ 
			\l(\frac{(s^2_1+s^2_2+s^2_3)\oa_3}{3s_3},0 \r) & \mbox{ if } t \in [t^f_{123},T]
		\end{cases}
	\end{equation*}
	and 
	\begin{equation*}
		\label{e46}
		\bar\x_i(t) = 
		\begin{cases}
			\bar \x_i(0) + (t\oa_i s_i, 0) & \mbox{ if } t \in [0,t^f_{123}) \\ 
			\bar\x_i(0) + \l(t^f_{123}\oa_is_i+ (t-t^f_{123})\frac{(s^2_1+s^2_2+s^2_3)\oa_3}{3s_3},0 \r)& \mbox{ if } t \in [t^f_{123},T],
		\end{cases}
	\end{equation*}
	for  $i=1, 2,3$.
	
	 Using the fact $\l\| \bar\x_3(t^f_{123})-\bar\x_2(t^f_{123})\r\| = L_2+L_3$ and  $\l\| \bar\x_2(t^f_{123})-\bar\x_1(t^f_{123})\r\| = L_1+L_2$ yields
	$$
	\begin{cases}
		(\oa_1s_1-\oa_2s_2)t^f_{123} = 2\Lm_{12} \\
		(\oa_2s_2-\oa_3s_3)t^f_{123} = 2\Lm_{23},
	\end{cases}
	$$
	which together with~\eqref{e30s} implies that 
	\begin{equation}
		\label{e47}
		\begin{cases}
			\frac{\Lm_{12}}{s_1^2-s_2^2} = \frac{\Lm_{23}}{s_2^2-s_3^2} \\ 
			t^f_{123} = \frac{2\Lm_{12}}{\oa_1s_1-\oa_2s_2} = \frac{2\Lm_{23}}{\oa_2s_2-\oa_3s_3}.
		\end{cases}
	\end{equation}
	
 The methodologies we have developed allow for the explicit calculation of the agents' final positions and the associated cost functional, as defined in~\eqref{e38}. It is important to note the critical role of the first equation in~\eqref{e47}. This equation acts as a condition to validate specific scenarios in our model. For instance, it can be straightforwardly confirmed that
\begin{enumerate}
	\item The condition $t^f_{12} < t^f_{23}$ holds true  if and only if $\frac{\Lm_{12}}{s_1^2-s_2^2} <\frac{\Lm_{23}}{s_2^2-s_3^2}$.
	\item Conversely $t^f_{12} > t^f_{23}$ occurs if and only if $\frac{\Lm_{12}}{s_1^2-s_2^2} >\frac{\Lm_{23}}{s_2^2-s_3^2}. $ 
\end{enumerate}	
These insights are particularly valuable as they provide a practical guideline to predict the sequence of contact among the agents. This reinforces the significance of the relationships established in equations~\eqref{e22} and~\eqref{e30s}, which are derived from the necessary optimality conditions presented in Theorem~\ref{Th3}. Our findings and the corresponding algorithmic approach are encapsulated in Algorithm~\ref{alg3} and further demonstrated through subsequent examples.

\begin{figure}[!ht]
		\centering
		\includegraphics[scale=0.4]{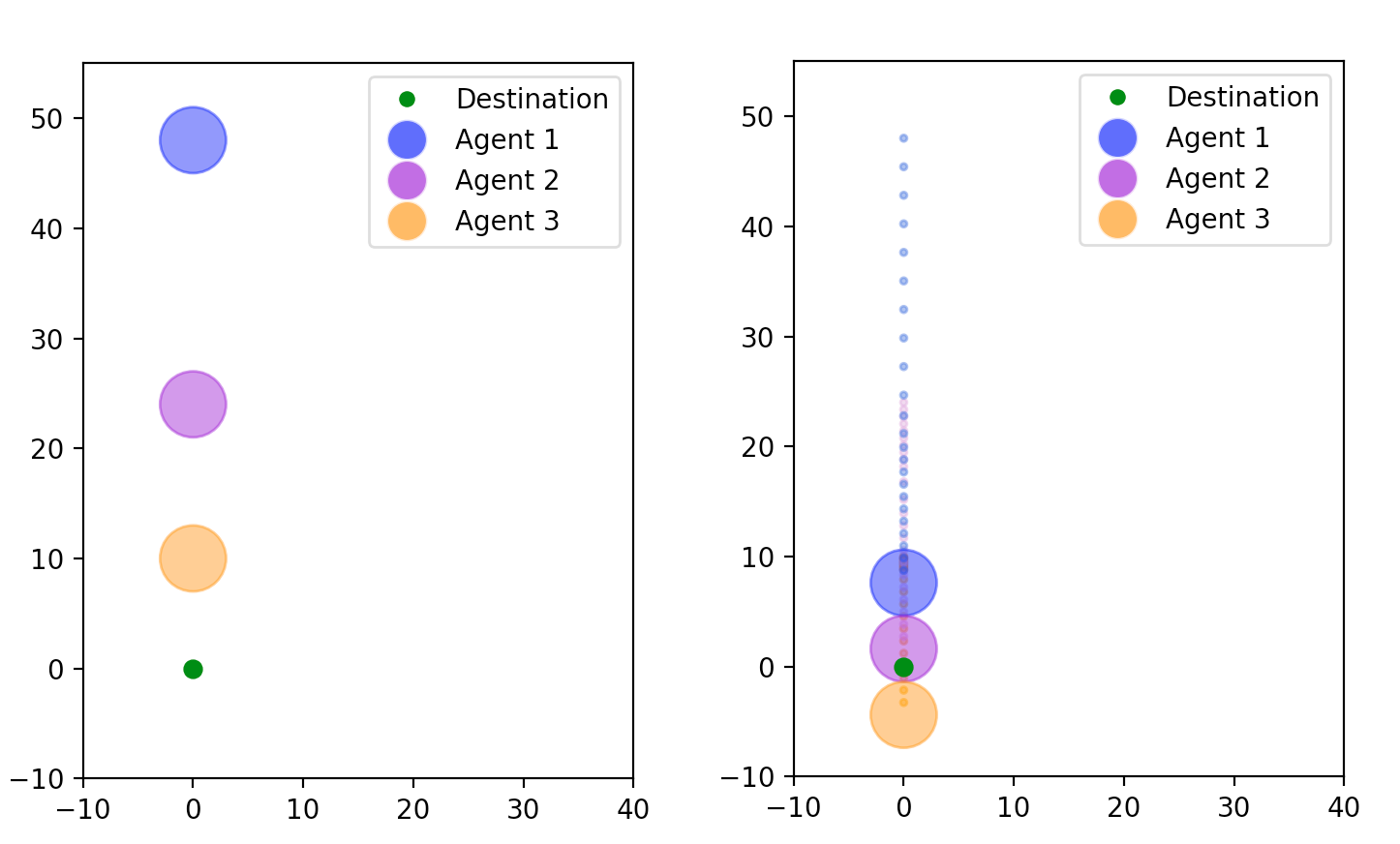}  
		\caption{Illustration for Example~\ref{eg53} at $t=0$ and $t=T$}\label{eg53fig}
	\end{figure}

	\begin{example}\label{eg53}
	We examine a controlled motion problem with the following configuration, as depicted  in Figure~\ref{eg53fig}:
$$
\begin{cases}
	\xd=(0,0), \; \x_1^0 = (0,48), \; \x_2^0= (0,24), \; \x_3^0=(0,10) \\
	T = 6, \; L_1=L_2=L_3=3.
\end{cases}
$$
The performances metrics for the agents, calculated under various conditions, are summarized in the table below.
		{\small \begin{center}
			\begin{tabular}{|c|c|c|c|c|c|c|}
				\hline
				$\tau$ &  $\oa_1$ &    $\oa_2$ &    $\oa_3$ & $t^f_{12}$ &  $t^f_{23}$  &       $J[\bar{x},\bar{a}]$ \\
				\hline
				1.0 &  1.306309 &  0.653154 &  0.272148 &  2.296547 &  2.796989 &   60.465 \\
				2.0 &  1.291812 &  0.645906 &  0.269128 &  2.322319 &  2.828377 &   84.660 \\
				3.0 &  1.277315 &  0.638658 &  0.266107 &  2.348676 &  2.860477 &  108.585 \\
				4.0 &  1.262819 &  0.631409 &  0.263087 &  2.375638 &  2.893314 &  132.240 \\
				5.0 &  1.248322 &  0.624161 &  0.260067 &  2.403226 &  2.926914 &  155.625 \\
				6.0 &  1.233825 &  0.616913 &  0.257047 &  2.431462 &  2.961303 &  178.740 \\
				7.0 &  1.219329 &  0.609664 &  0.254027 &  2.460370 &  2.996510 &  201.585 \\
				8.0 &  1.204832 &  0.602416 &  0.251007 &  2.489973 &  3.032564 &  224.160 \\
				9.0 &  1.190336 &  0.595168 &  0.247987 &  2.520298 &  3.069497 &  246.465 \\
				10.0 &  1.175839 &  0.587919 &  0.244966 &  2.551370 &  3.107340 &  268.500 \\
				\hline
			\end{tabular}
		\end{center}}
	\end{example}

\begin{remark}
	\label{3to2}
	When two of the three agents, such as $\bar\x_1$ and $\bar\x_2$, are initially in contact, they can be collectively represented by a single agent, denoted as 
		$\bar\x_{12}:= \frac{1}{2}(\bar\x_1+\bar\x_2)$.
		This simplification effectively models them as an elastic with center $\bar\x_{12}$ and radius $L_{12} = L_1+L_2$. Consequently, the three-agent model can be reduced to a two-agent scenario. In such cases, agent 2 assumes the role of the leader, with agent 1 mimicking its actions. Although models with more than three agents introduce greater complexity, they can be similarly streamline by appropriately grouping the agents into three categories. 
\end{remark}

\section{Concluding Remarks and Future Research}
This paper has delved into various optimal control problems related to crowd motion models, deriving necessary optimality conditions that facilitate systematic problem-solving. These conditions, combined with constructive algorithms, enable us to address dynamic optimization problems involving single and multiple agents in corridor settings. While some problems, particularly those involving a singe agent and an obstacle, are not entirely resolved, our optimality conditions still offer valuable insights for devising optimal control strategies. The models and frameworks we have developed hold potential for applications in real-world scenarios, paving the way for future research in this domain.

\section{Algorithms}
\setcounter{equation}{0}
This section presents two algorithms used in our paper.

	\begin{breakablealgorithm}\label{alg2}
		\caption{Optimal control for the crowd motion models with two agents in a corridor }
		\begin{algorithmic}[1] 
			\Procedure {optimal a}{$T$, $\x^0$, $\xd$, $L_1$, $L_2$, $\tau $}
			\State Compute $s_i =  \frac{\| \xd - \bar\x_i(0)\|}{T}$
			\If {$\bar\x_{21}(0)-\bar\x_{11}(0)=L_1+L_2$}
			\State minimize  $J[\bar\x, \oa]$ in~\eqref{e32}
			\Else
			\State minimize  $J[\bar\x, \oa]$ in~\eqref{e33} (there is no contact)
			
			\State compute $\eta_{12}(t^f_{12})$ in terms of $\oa_2$ using~\eqref{e26}
			\State compute $\Lm_{12}$ and $t^f_{12}$ in terms of $\oa_2$ using~\eqref{e27}
			\State compute $\x_1(t), \x_2(t)$ in terms of $\oa_2$ using~\eqref{e26a}
			\State minimize  $J[\bar\x, \oa]$ in~\eqref{e28}
			\State compare the minimum cost in~\eqref{e30} and~\eqref{e33}
			\EndIf 
			\State compute $\oa_1 = \frac{s_1}{s_2}\oa_2$
			\EndProcedure
		\end{algorithmic}
	\end{breakablealgorithm}
	
	\begin{breakablealgorithm}\label{alg3}
		\caption{Optimal control for the crowd motion models with three agents in a corridor }
		\begin{algorithmic}[1] 
			\Procedure {optimal a}{$T$, $\x^0$, $\xd$, $L_1$, $L_2$, $L_3$, $\tau $}
			\State compute $s_i =  \frac{\| \xd - \bar\x_i(0)\|}{T}$
	
			\State compute $\Lm= \frac{\Lm_{12}}{s_1^2-s_2^2} -\frac{\Lm_{23}}{s_2^2-s_3^2}$
			\If{$\Lm <0$}
			\State $t^f_{12} < t^f_{123}=t^f_{23}$
			\State compute $\eta_{12}(t^f_{12}), \eta_{12}(t^f_{123}),$ and $\eta_{23}(t^f_{123})$ using~\eqref{e34s} and~\eqref{e35s}
			\State compute $t^f_{123}$ using~\eqref{e31a}
			\State compute $\bar\x_1(T), \bar\x_2(T),$ and $\bar\x_3(T)$ using~\eqref{e37a}--\eqref{e37c}
			\State minimize $J[\bar\x,\oa]$ in~\eqref{e38}
			\EndIf
			\If{$\Lm>0$}
			\State $t^f_{23}<t^f_{123}=t^f_{12}$
			\State compute $\eta_{23}(t^f_{23}), \eta_{23}(t^f_{123}),$ and $\eta_{12}(t^f_{123})$ using~\eqref{e39} and~\eqref{e40}
			\State compute $t^f_{123}$ using~\eqref{e31b}
			\State compute $\bar\x_1(T), \bar\x_2(T),$ and $\bar\x_3(T)$ using~\eqref{e42a}--\eqref{e42c}
			\State minimize $J[\bar\x,\oa]$ in~\eqref{e38}
			\EndIf
			\If{$\Lm=0$} 
			\State compute $\eta_{12}(t^f_{123})$ and $\eta_{23}(t^f_{123})$ using~\eqref{e43} and~\eqref{e44}
			\State minimize  $J[\bar\x,\oa]$ in~\eqref{e38}
			\EndIf
			\EndProcedure
		\end{algorithmic}
	\end{breakablealgorithm}

\noindent \textbf{Acknowledgements.}
	The authors are indebted to Professor Boris Mordukhovich for his helpful remarks and discussions on the original presentation. The authors would also like to express their sincere appreciation to Kangmin Cho, Jinwoo Choi, Sinae Hong, Abhishek Kafle, Hansol Lim, Biniam Markos, and Jiung Seo for their suggestions on the research project. Research of Tan~H.~Cao was supported by the National Research Foundation of Korea grant funded by the Korea Government (MIST) NRF-2020R1F1A1A01071015.

\end{document}